\documentclass[a4paper,11pt]{amsart}
\usepackage[top=2.2cm,left=2.8cm,right=2.8cm,bottom=2.8cm]{geometry}
\usepackage[foot]{amsaddr}

\usepackage{hyperref}
\hypersetup{
	colorlinks,
	citecolor=green,
	filecolor=black,
	linkcolor=blue,
	urlcolor=blue
}
\hypersetup{linktocpage}

\usepackage{graphicx}

\usepackage{amsbsy}
\usepackage{latexsym}
\usepackage{amsfonts}
\usepackage{amssymb}
\usepackage{cleveref}
\usepackage{upgreek}
\usepackage[usenames]{color}
\usepackage{amsmath,amsthm}
\usepackage{enumerate}
\usepackage{pdfpages}

\usepackage{tikz}
\usetikzlibrary{arrows.meta,calc,matrix}
\usetikzlibrary{decorations.pathreplacing}

\usepackage{stmaryrd}
\usepackage{dsfont}
\usepackage{bbm}

\DeclareMathOperator{\ran}{ran}

\DeclareMathOperator{\rank}{rank}

\DeclareMathOperator{\Span}{span}

\newtheorem{theorem}{Theorem}[section]
\newtheorem{lemma}[theorem]{Lemma}
\newtheorem{prop}[theorem]{Proposition}

\theoremstyle{definition}

\newtheorem{example}[theorem]{Example}

\theoremstyle{remark}
\newtheorem{remark}[theorem]{Remark}

\numberwithin{equation}{section}

\theoremstyle{plain}

\newcommand{\Chi}{\raise .3ex
	\hbox{\large $\chi$}}

\newcommand{\R}{\mathbb{R}}

\newcommand{\N}{\mathbb{N}}

\newcommand{\eps}{\varepsilon}
\usepackage[section]{algorithm}
\usepackage{algorithmicx}
\usepackage{algpseudocode}

\usepackage[maxbibnames=50,style=alphabetic]{biblatex}
\title[Accurate rank-controlled step-truncation methods]{Accurate rank-controlled step-truncation methods for matrix differential equations}

\author{Matthieu Dolbeault${}^1$}
\address{${}^1$ CNRS and Laboratoire de Mathématiques Jean Leray, Nantes Université, France}
\email{matthieu.dolbeault@univ-nantes.fr}
\author{Polina Sachsenmaier${}^2$}
\address{${}^2$ Institut für Geometrie und Praktische Mathematik, RWTH Aachen University, Germany}
\email{sachsenmaier@igpm.rwth-aachen.de}
\author{Fabio Zoccolan${}^3$}
\address{${}^3$ Institut de Mathématiques, École Polytechnique Fédérale de Lausanne, Switzerland}
\email{fabio.zoccolan@epfl.ch}

\begin{document}

\begin{abstract}
A class of explicit time stepping schemes is introduced for the integration of matrix differential equations in low-rank format. These schemes combine explicit Runge--Kutta time stepping with rank-adaptive recompressions. For suitably chosen truncation thresholds, they converge at arbitrary order, while having close to optimal ranks, without any prior assumptions on the best approximation ranks. The ranks of the scheme also necessarily depend on the time step, as shown by  counterexamples. A numerical experiment illustrates the advantages and performance of the proposed scheme.


\end{abstract}

\maketitle

\section{Introduction}

High dimensional systems of ordinary differential equations (ODEs) arise in many practical applications, ranging from interacting particle systems to parametric and stochastic problems, as well as for the training of huge machine learning architectures. They also constitute a first step towards
the numerical study of partial differential equations.
Numerical methods for such differential equations are often limited by their computational cost and storage requirements. This calls for model order reduction techniques, which leverage a lower-dimensional structure in the unknown solution.

In this work, we consider ODEs of the form
\begin{equation}
\partial_t X(t) = F(t,X(t)), \qquad X(0) = X_0,
\label{eq: ODE}
\end{equation}
where the unknown $X(t)\in \R^{d_1\times d_2}$ is a matrix of large dimensions $d_1\geq d_2\gg1$ evolving over the time interval $t\in [0,T]$, and the right-hand side $F$ is a smooth function from $[0,T]\times \R^{d_1\times d_2}$ to $\R^{d_1\times d_2}$. A direct numerical integration of \eqref{eq: ODE} requires $\mathcal{O}(d_1 d_2)$ operations per time step. In many cases of interest, however, the singular values of the solution $X(t)$ decay rapidly, so that it can be efficiently approximated by matrices of rank $r \ll d_2$. Methods operating only on the corresponding low-rank factors then require $\mathcal O(d_1 r^2)$ operations per time step, that is, their cost scales linearly in the matrix dimensions, rather than quadratically. In contrast to projector-based model order reduction, see for example \cite{quarteroni2015reduced} and \cite{rowley2017model}, where approximations are sought in fixed subspaces constructed in an expensive offline phase, the low-rank approximations considered here are updated in each step. 

The accuracy and efficiency of such methods depend crucially on the ranks of the approximation. This involves a balance between two competing requirements: the ranks need to be sufficiently large to keep approximation errors small, but should at the same time not be much larger than necessary in view of the computational costs. Since the singular values of $X(t)$ are time-dependent, an appropriate rank is in general not known in advance and may vary over the time interval, making an adaptation of the ranks in the course of the time integration necessary. 

Existing low-rank time integration methods rely on different mechanisms for controlling ranks. On the one hand, methods such as dynamical low-rank approximation \cite{koch2007dynamical} keep the ranks fixed, but their projection errors cannot be controlled in general without further assumptions. On the other hand, step-truncation methods \cite{ST_2022} are able to approximate standard time-stepping schemes up to a prescribed accuracy, but the resulting ranks may also be substantially larger than necessary. Guaranteed error bounds combined with bounds on the ranks in terms of best approximation ranks have so far only been established for methods based on space-time variational formulations for parabolic equations \cite{BachmayrFaldum24} and for iterative thresholding schemes for linear Schr\"odinger equations \cite{BDS25}, \cite{BJSV26}.

The main result of this work is a modified step-truncation method based on explicit Runge-Kutta schemes, for which we establish both convergence and control of the ranks. More precisely, the method converges with order $q$, for any $q \in \mathbb N$, in the time step $h$ and at the same time, the ranks of all computed quantities remain comparable to the minimal ranks of approximations of the exact solution with the same accuracy. To the best of our knowledge, no result of this type has previously been available for step-truncation methods.

For each fixed $t$, the best low-rank approximation of $X(t)$ is given by truncation of its singular value decomposition (SVD)
\[
X=\sum_{k=1}^{d_2}\sigma_k u_k v_k^\top,
\]
for a non-increasing sequence of singular values $\sigma_k\geq 0$, and orthonormal families of singular vectors $u_k\in \R^{d_1}$ and $v_k\in \R^{d_2}$. Given the SVD framework, it is natural to measure the accuracy of the approximation in terms of the Schatten $p$-norm
\begin{equation}
\label{eq:schatten}
\|X\|_p=\bigg(\sum_{k=1}^{d_2}\sigma_k^p\bigg)^{1/p}, \qquad p\in [1,\infty],
\end{equation}
where $\|X\|_\infty = \sigma_1$ is the standard spectral norm and $\|X\|_F := \|X\|_2$ the Frobenius norm. For $\delta>0$, we define the \textit{recompression operator} $\mathcal R_\delta: \mathbb R^{d_1 \times d_2} \to \mathbb R^{d_1 \times d_2}$ as
\begin{equation}\label{eq: R_delta}
    \mathcal R_\delta(X) = \sum_{k=1}^{r_\delta} \sigma_k u_kv_k^\top, \quad \text{where} \quad r_\delta = \min \Big\{ r \geq 0 \, : \, \sum_{k > r} \sigma_k^p \leq \delta^p \Big\}.
\end{equation}
Note that for $p = \infty$ the condition reads $\sigma_{r+1} \leq \delta$. As a consequence of the Eckart-Young-Mirsky theorem \cite{mirsky1960symmetric}, $\mathcal R_\delta(X)$ has minimal rank among all approximations of $X$ with error at most $\delta$, that is, 
\begin{equation*}
    \min_{\|Y-X\|_p \leq \delta}\rank(Y) = r_\delta.
\end{equation*}
Even though the SVD and thus $\mathcal R_\delta$ are not unique when $\sigma_{r_\delta} = \sigma_{r_\delta+1}$, our results hold for any choice. For each given accuracy $\delta$, the ranks $r_\delta$ are thus the natural benchmark for the ranks of any numerical approximation. Note that both $r_\delta$ and $\mathcal R_\delta$ depend on the choice of~$p$; in what follows, the operator $\mathcal R_\delta$ is thus always understood with respect to the norm in which errors are measured.

Our aim is to construct an approximation $X_n \approx X(t_n)$ whose rank is comparable, up to a moderate constant, to the rank of the best approximation of $X(t_n)$ with the same accuracy. In classical step-truncation methods, the truncation tolerances are chosen proportional to the local error of the time-stepping scheme and thus decrease with $h$, and no comparison of the resulting ranks with best approximation ranks is available. We show that this difficulty can be overcome by using two recompression procedures with tolerances proportional to the current time-stepping error, selected and placed such that the global error does not grow exponentially with the number of time steps.

\subsection{Main result}\label{sec: main res}

Let $h>0$ denote the time step size. For convenience, we assume that $1/h\in \N$ throughout the article. We discretize the time interval $[0,T]$ as a sequence $(t_n)_{0\leq n\leq N}$, with $t_n=nh$ and $N=\lfloor T/h\rfloor$. We also fix $p\in [1,\infty]$ and let $\|\cdot\|=\|\cdot\|_p$.

For an arbitrary $q\in \N$, we construct a scheme $X_n$ approximating the solution $X(t_n)$ at order $q$, that is, with an error of order $O(h^q)$. The proposed scheme has \textit{near-optimal ranks}, in the sense that the rank of $X_n$ is bounded up to a multiplicative factor by the minimal rank of any equally accurate approximation of $X(t_n)$. 
Our strategy consists in embedding low-rank recompressions into the iterations of an explicit time integrator $\Psi(t,X)$, satisfying
\[
\Psi(t,X(t))\approx \frac{X(t+h)-X(t)}{h}= \frac1h \int_t^{t+h} F(s,X(s)) \, ds,\qquad t\in [0,T].
\]

We make the following assumptions: first, as for any numerical method of global order $q$, the time integrator should have local accuracy of order $q+1$, that is,
\begin{equation}
\|X(t+h)-X(t)-h\Psi(t,X(t))\|\leq Ch^{q+1},\qquad t\in[0,T-h].
\label{assumption-local-Psi}\tag{A1}
\end{equation}
Throughout the paper, the letter $C$ will refer to the positive constant appearing in assumption~\eqref{assumption-local-Psi}.
Second, to propagate error from one time step to the next, $\Psi$ is assumed to be Lipschitz with respect to the second argument with some constant $\Lambda>0$ 
\begin{equation}
\|\Psi(t,X)-\Psi(t,Y)\| \leq \Lambda\,\|X-Y\|,\qquad t\in[0,T],\quad X,Y\in \R^{d_1\times d_2}.
\label{assumption-lipschitz-Psi}\tag{A2}
\end{equation}
Finally, applying $\Psi$ should almost preserve the low-rank structure, in the sense that there exist $R_\Psi\in \N_0$ and $R_0: [0,T] \times (0, \infty) \rightarrow \mathbb N_0$ such that for all $t$, $X$ and $\gamma > 0$ there is a $Z \in \mathbb R^{d_1 \times d_2}$ with
\begin{equation}
\| Z - \Psi(t,X) \| \leq \gamma, \qquad \rank (Z) \leq R_\Psi \rank (X) + R_0(t, \gamma).
\label{assumption-rank-Psi}\tag{A3}
\end{equation}
We also define $R_0(t,0) = \sup_{\gamma > 0}R_0(t, \gamma)$. Note that if the right-hand side $F$ does not incorporate a source term, one can set $R_0 = 0$.  For computational purposes, we moreover assume that an SVD of $\Psi(t,X)$ can be obtained from an SVD of $X$ in $O \big(d_1\rank(X)^2 \big)$ operations, and refer to Remark~\ref{rk:ranks-depend-on-h}.

\begin{theorem}
\label{thm:main}
Under assumptions \eqref{assumption-local-Psi}--\eqref{assumption-rank-Psi},  define for $0 \leq n \leq N$ the target accuracy and optimal approximations ranks as
\begin{equation}
    \xi_n = Ch^q\, 3^{j_n} e^{3(\Lambda+1)t_n},
    \qquad
    \bar r_n = \min_{\|Y-X(t_n)\|_p \leq \xi_n} \rank(Y),
    \label{eq: thm 1 opt ranks}
\end{equation}
where $j_n = \max\{\lceil t_n\rceil - 1, 0\}$ is the number of positive integers smaller than~$t_n$. Then the approximations produced by Algorithm~\ref{alg:step_trunc_reset} satisfy the following bounds.
\begin{enumerate}[(i)]
    \item \emph{Accuracy:} For $0 \leq n \leq N$, it holds that
    \begin{equation*}
        \|X_n - X(t_n)\|_p \leq \xi_n.
    \end{equation*}
    
    \item \emph{Ranks:} For $0 \leq n \leq N-1$, one has the rank bound
    \[
        \rank(X_{n+1}) \leq \frac{1}{h}\Big(1 + \max_{n' \leq k \leq n} \big( e^{-3^{-p}h^{p-1}(t_n-t_k)}\, (1+R_\Psi)\,\bar r_k + R_0(t_k,\xi_k) \big) \Big) + R_0(t_n,0),
    \]
    where $n' = \lfloor t_n\rfloor/h$ is the index of the last integer time.
    \item \emph{Ranks after a reset:} For integer times $t_n \in \{0,1,2,\dots\}$, we have the sharper rank bound
    \[
        \rank(X_{n+1}) \leq (1+R_\Psi)\,\bar r_n + R_0(t_n,0), \qquad 0 \leq n \leq N-1.
    \]
\end{enumerate}
\end{theorem}

Before comparing our result with the literature, a few explanatory remarks are in order. First, the assumptions can be derived from properties of $F$, see Section~\ref{sec: RK ODEs}. The constants $C$, $\Lambda$ and $R_\Psi$ as well as the function $R_0$ can be related to the smoothness, Lipschitz stability, and ranks of the right-hand side $F$ of equation~\eqref{eq: ODE}. In particular, since $R_0$ accounts for source terms, it vanishes for homogeneous problems. Assumption \eqref{assumption-rank-Psi} then reduces to $\rank(\Psi(t,X)) \leq R_\Psi \rank(X)$. For right-hand sides of the form $F(t,X) = F_0(t,X) + G(t)$ it is governed by the best approximation ranks of $G$.

Second, for $\xi_0 = Ch^q$, the target error $\xi_n$ has the desired order $O(h^q)$, with the same constant $C$ as in \eqref{assumption-local-Psi}, and the usual exponential dependence in the time $t_n$, as discussed in Remark~\ref{rem: exp T dependence}. The target rank $\bar r_n$ corresponds to the minimal rank of any approximation of $X(t_n)$ with error at most~$\xi_n$. As we do not know the exact solution $X(t_n)$, we cannot expect our scheme $X_n$ to exactly reach the optimal rank $\bar r_n$. Nevertheless, the rank of $X_{n+1}$ is bounded by a multiple of the $\bar r_k$ at previous times $t_k$ of the current unit interval, plus terms accounting for the source, justifying the claimed near-optimality. In particular, we have
\[
\rank(X_{n+1})\leq \frac{1+R_\Psi}{h} \Big(1+\max_{k\leq n} \, \bar r_k \Big) + \frac 1h \max_{k \leq n} R_0(t_k, \xi_k) + R_0(t_n,0).
\]

Third, the terms in the rank bound admit a natural interpretation: past ranks are damped exponentially, at a rate that depends on the choice of the norm, see Remark~\ref{rk:choice-p}, and the dependence is limited to times between $\lfloor t_n\rfloor$ and $t_n$. The factor $R_\Psi$ reflects the use of assumption~\eqref{assumption-rank-Psi}. The term $R_0(t_k, \xi_k)$ takes into account the part of the source, which is relevant at the target accuracy and enters in every step. In contrast, the exact rank of the latest source contribution, that is, $R_0(t_n,0)$, appears only once, since $X_{n+1}$ contains the last increment without recompression, see Algorithm \ref{alg:step_trunc_reset}. The additional ones prevent issues in the trivial cases $R_\Psi=0$ or $\bar r_k=0$. Note that if $\bar r_n\geq 1$, then the additive $1$ in the rank bound may be absorbed into the maximum at the cost of a factor at most 2.

\begin{remark}
\label{rk:ranks-depend-on-h}
In Theorem~\ref{thm:main}, the ranks are proven to be near-optimal up to a factor $1/h$. This dependence in $h$ is in fact necessary, as shown in Proposition \ref{prop: h sharpness}. Moreover, since the computational cost of the algorithm scales as $\rank(X_n)^2$, this factor only corresponds to an increase of the computational cost from $O(T/h)$ to $O(T/h^3)$. This is equivalent to a reduction of the convergence rate, in terms of scaling between error and computational effort. In contrast, the reduction to a linear dependence in the matrix dimensions is crucial for large values of $d_1$ and $d_2$. It will play an even more important role in the extension of the present approach to higher order tensor formats, in the spirit of \cite{BJSV26}, by avoiding the \emph{curse of dimensionality}, that is, computational and storage costs scaling exponentially in the tensor order.
\end{remark}

\begin{remark}
\label{rk:choice-p}
As can be noticed in the statement of Theorem~\ref{thm:main}, the exponential factor multiplying $\bar r_k$, indicating the vanishing importance of past best approximation ranks, depends on the choice of norm. For $p=1$, corresponding to the trace norm $\|X\|_1=\sum\sigma_k$, this exponential decay is proportional to the elapsed time. In practice, if the ranks of the solution decrease, we expect the ranks of the scheme to decrease accordingly. The case of the Frobenius norm $p=2$ was investigated in \cite{BDS25} for time-dependent linear Schr\"odinger equations, since it is the natural norm for a discretization of the $L^2$ norm on multivariate domains; with this choice the analysis only provides a slower decay. Finally, in the spectral norm $p=\infty$, with the convention $h^{p-1}=0$ hence $\exp(-3^{-p}h^{p-1}(t_n-t_k))=1$, the result still holds but without exponential decay; in this case we rely on the fact that the ranks are reset to a near-optimal value every $1/h$ steps. 
\end{remark}

\subsection{Novelty and relation to prior works} 

The basic ingredients of the proposed scheme, explicit Runge-Kutta methods and recompression by truncated SVD, are standard. The main novelty lies in the placement of the recompressions and the particular choice of the corresponding tolerances, which together ensure convergence of order $q$ as well as rank bounds in terms of best approximation ranks.

Bounds of this type have previously been obtained for methods based on space-time variational formulations, in which an approximation of the entire evolution is refined jointly in the spatial and temporal variables. Low-rank methods of this kind were considered in \cite{andreev2015multilevel}, \cite{boiveau2019low} as well as \cite{dupuy2024space} and near-optimal rank bounds were established for high-dimensional parabolic problems in \cite{BachmayrFaldum24}.

Well-posed space-time formulations, however, are not available for all problems of interest. The iterative thresholding scheme of \cite{BDS25}, extended to higher-order tensors in \cite{BJSV26}, does not require such a formulation: on each subinterval, approximations are computed by fixed-point iterations for implicit collocation methods combined with soft-thresholding of the iterates, which yields quasi-optimal rank bounds for linear Schr\"odinger-type problems. This scheme builds on adaptive iterative methods for high-dimensional operator equations, in which recompressions with tolerances adapted to the current error lead to near-optimal ranks, see \cite{BachmayrDahmen15} and \cite{bachmayr2017iterative}.

For methods that are not of this type, to the best of our knowledge, either convergence or rank control can be ensured, but not both. Dynamical low-rank approximation (DLRA) \cite{koch2007dynamical} belongs to the second class: the right-hand side is projected onto the tangent space of the manifold of matrices of an initially fixed rank at the current approximation, and the error bounds rely on the assumption that the resulting projection error is small, which is not guaranteed in general. This observation appears in \cite[Remark 7.2]{bachmayr2023low}, as well as Example~1 of \cite{lam2025randomized}, which adapts the DLRA framework to arbitrary order Runge-Kutta schemes. The projector-splitting \cite{lubich2014projector} and basis update \& Galerkin (BUG) \cite{ceruti2022unconventional} integrators are robust with respect to small singular values \cite{kieri2016discretized}. Moreover, the BUG integrator can be extended to a rank-adaptive version \cite{ceruti2022rank}, \cite{ceruti2024parallel}. Further rank-adaptive and high-order variants were proposed in \cite{hauck2023predictor}, \cite{charous2023dynamically}, \cite{charous2024stable}, \cite{nobile2026high}, \cite{carrel2025interpolatory} and \cite{kieri2019projection}. We also refer to \cite{seguin2024low} for the relation to retractions and to \cite{einkemmer25} for a recent review. The error analysis of all these methods, however, still depends on assumptions on the tangent space projection error.

In contrast, step-truncation methods use low-rank representations of the right-hand side: a step of a standard time-stepping scheme, such as an explicit Runge-Kutta method, is performed in low-rank format, and the result is subsequently truncated, see \cite{ST_2022}, \cite{lam2025randomized} and \cite{guo2022low}. Implicit schemes of this type were considered in \cite{ST_2023}, \cite{li2024high} and \cite{nakao2025reduced}, where the latter can be viewed as a hybrid between DLRA and step-truncation. Methods using the right-hand side for basis enrichment are presented in \cite{appelo2025robust} and \cite{yang2020time}. In order to retain the order of the underlying scheme, the truncation tolerances are chosen in terms of its local error \cite{ST_2023} and therefore decrease with the time step. While convergence can thus be ensured, no bounds on the resulting ranks are available. In the present work, we show that additional recompressions lead to rank bounds in terms of best approximation ranks for such methods. In comparison to \cite{BDS25} and \cite{BJSV26}, the resulting scheme requires neither fixed-point iterations nor implicit time stepping, and applies to general ODEs of the form \eqref{eq: ODE} satisfying \eqref{assumption-local-Psi}--\eqref{assumption-rank-Psi}.

\subsection{Outline}

In Section~\ref{sec: RK ODEs}, we recall properties of explicit Runge--Kutta methods, relate assumptions~\eqref{assumption-local-Psi}--\eqref{assumption-rank-Psi} to properties of $F$, and discuss classical step-truncation methods. In Section~\ref{sec: step-trunc add recom}, we introduce the proposed scheme in Algorithm~\ref{alg:step_trunc_reset} and prove Theorem~\ref{thm:main}. The optimality of the rank bounds is discussed in Section~\ref{sec: optimality of the bounds}, and a numerical experiment illustrating the advantage of the proposed scheme over the classical rank-adaptive BUG approach is presented in Section~\ref{sec: num exp}.

\section{Step-truncation methods for ordinary differential equations}
\label{sec: RK ODEs}

The aim of this work is to develop an efficient and easy to implement explicit time-stepping scheme in the low-rank format with near-optimal rank guarantees. More specifically, we propose and analyze a modified step-truncation method based on explicit Runge-Kutta methods with recompressions placed so that an exponential growth of the global error with respect to the number of time steps is avoided, while simultaneously maintaining rank control on all computed quantities.

To this end, we assume in the following that the right-hand side $F$ of equation \eqref{eq: ODE} satisfies the following: first of all, the right-hand side is assumed to be smooth enough, in particular
\begin{equation}
\label{eq: F smooth}\tag{A4}
F\in C^q([0,T]\times \R^{d_1\times d_2},\R^{d_1\times d_2}),
\end{equation}
Lipschitz stable with constant $L>0$, that is,
\begin{equation}
\label{eq: F lip}\tag{A5}
\|F(t,X)-F(t,Y)\|\leq L\|X-Y\|,\qquad t\in [0,T],\quad X,Y\in \R^{d_1\times d_2},
\end{equation}
and satisfies a rank bound up to accuracy $\gamma$ of the form
\begin{equation}
\label{eq: rank F bound}\tag{A6}
\rank \Big(\mathcal R_{\gamma}\big(F(t,X)\big) \Big)\leq R_F\rank(X) + R_G(t, \gamma),\qquad \gamma >0, \quad t\in [0,T],\quad X\in \R^{d_1\times d_2},
\end{equation}
for some $R_F \in \mathbb N_0$ and $R_G: [0,T] \times (0,\infty) \rightarrow \mathbb N_0$.
For nonlinear problems, one can restrict these conditions to a bounded neighbourhood of the exact trajectory $\big\{ X(t)  \, : \,  t \in [0,T] \big\}$. Since all arguments at which $F$ is evaluated remain close to the trajectory for sufficiently small $h$, only a local Lipschitz property is needed in \eqref{eq: F lip}.

While the first two assumptions are classical, see for instance \cite{hairer1993solving}, the third one is essential in order to be able to establish rank control on the produced iterates of the algorithm. In particular, these assumptions include affine right-hand sides of the form
\begin{equation*}
    F(t,X) = \sum_{i=1}^m A_i(t) X B_i(t) + G(t),
\end{equation*}
where the matrix-valued functions $A_i$, $B_i$ and $G$ are of class $C^q$. For such right-hand sides \eqref{eq: F lip} holds with $L = \sum_i \sup_t \| A_i(t) \|_\infty \|B_i(t)\|_\infty$ and \eqref{eq: rank F bound} holds with $R_F = m$ and $R_G(t,\gamma) = \rank \big(\mathcal R_\gamma(G(t)) \big)$.

\subsection{Runge-Kutta methods for ordinary differential equations}

We start by recalling properties of explicit Runge-Kutta methods, for which we refer to \cite{iserles}. The general strategy to approximate the solution of \eqref{eq: ODE} can be written in the additive one-step form
\begin{equation*}
    Z_{n+1} = Z_n + h \Psi(t_n, Z_n),
\end{equation*}
where $Z_n \approx X(t_n)$ for $1 \leq n \leq N$ is the numerical approximation of $X$ at time $t_n$. Given a Butcher tableau $(a_{i,j},b_j,c_i)$ with $a_{i,j}=0$ for $j \geq i$ and $c_1 = 0$, the Runge-Kutta stage matrices for an explicit Runge-Kutta method are defined as
    \begin{equation} \label{eq:RK_stages}
        \Psi^i(t,X) = F\bigg(t+ c_ih, X + h\sum_{j < i} a_{i,j} \Psi^j(t,X)\bigg), \qquad 1\leq i\leq s,
    \end{equation}
    where $s \in \mathbb{N}$ is the number of stages.
    Then, for all $h$, the operator $\Psi$ is given by
    \begin{equation*}
        \Psi(t, X) = \sum_{j=1}^s b_j \Psi^j.
    \end{equation*}
In the examples we will mostly use the explicit Euler scheme, the midpoint rule and RK4, which are of order 1, 2 and 4, respectively.

A Runge-Kutta method is said to be of order~$q$, if for every sufficiently smooth right-hand side $F$ (see Assumption \eqref{eq: F smooth}), the local error satisfies
\begin{equation*}
    \| Z_{n+1} - \tilde Z(t_{n+1}) \| = \mathcal O(h^{q+1}),
\end{equation*}
where $\tilde Z(t_{n+1})$ denotes the solution of \eqref{eq: ODE} with initial condition $\tilde Z(t_n) = Z_n$. Under Assumption \eqref{eq: F lip}, this then implies the global error bound
\begin{equation*}
    \max_{0 \leq n \leq N} \| X(t_n) - Z_n \| = \mathcal O(h^q).
\end{equation*}

Moreover, the above assumptions on the right-hand side $F$ imply those on the numerical scheme $\Psi$ stated in the introduction.
\begin{prop}
\label{prop: assumptions on F ot Psi}
Let $\Psi$ be an explicit $s$-stage Runge-Kutta method of order $q$. Then assumptions \eqref{assumption-local-Psi}, \eqref{assumption-lipschitz-Psi} and~\eqref{assumption-rank-Psi} follow from assumptions~\eqref{eq: F smooth}, \eqref{eq: F lip} and~\eqref{eq: rank F bound} with the Lipschitz constant
\begin{equation*}
        \Lambda = \sum_{i=1}^s|b_i| \Lambda_i, \qquad \Lambda_i = L + hL \sum_{j<i} |a_{i,j}|\Lambda_j,
\end{equation*} 
the rank bound of the homogeneous part
\begin{equation*}
    R_{\Psi} = \sum_{i,\ b_i \neq 0} R^i, \qquad R^i = R_F +  R_F\sum_{j < i, \ a_{i,j} \neq 0}R^j,
\end{equation*}
as well as the source bound
\begin{equation*}
    R_0(t, \gamma) = R^\prime \max_{1 \leq i \leq s} R_G \Big( t+c_ih, \frac{L}{\Lambda}\gamma \Big), \quad R^\prime = \sum_{i, \ b_i \neq 0} \rho_i , \quad \rho_i = 1 + R_F \sum_{j<i, \ a_{i,j} \neq 0} \rho_j.
\end{equation*}
for $1 \leq i \leq s$. Here, we use the convention that empty sums are equal to zero, such that $\Lambda_1 = L$ and $R^1 = R_F$.
\end{prop}

\begin{proof}
Using \eqref{eq: F smooth}, the exact solution satisfies $X \in C^{q+1}([0,T])$, so~\eqref{assumption-local-Psi} follows from \cite[Section II.3 Theorem 3.1]{hairer1993solving} using a classical local error estimate argument. The constant then depends only on the coefficients of the Runge-Kutta method and on the bounds for the partial derivatives of the right-hand side $F$ up to order $q$ on a bounded set containing the trajectory and all stage arguments.

Let further $X,Y\in\R^{d_1\times d_2}$ and $t \in [0,T]$. Using~\eqref{eq: F lip} and~\eqref{eq:RK_stages}, we obtain
\[
\|\Psi^i(t,X)-\Psi^i(t,Y)\| \leq L\Big(\|X-Y\| + h\sum_{j<i}|a_{i,j}|\,\|\Psi^j(t,X)-\Psi^j(t,Y)\|\Big),
\]
so that $\|\Psi^i(t,X)-\Psi^i(t,Y)\|\leq \Lambda_i\|X-Y\|$ by induction on $i$. Assumption \eqref{assumption-lipschitz-Psi} then follows from the triangle inequality, see also \cite[Lemma~1]{Lawson}.

Let $\gamma > 0$, and $Z \in \mathbb R^{d_1 \times d_2}$. We then define the approximate stages recursively as
\begin{equation*}
    W = \sum_i b_i W_i, \quad W_i = \mathcal R_{\gamma^\prime} \Big( F \big(t+c_ih, Z+h\sum_{j < i} a_{i,j} W_j \big) \big), \qquad \gamma^\prime = \frac{L}{\Lambda} \gamma.
\end{equation*}
It then holds that
\begin{equation*}
    \| \Psi^i - W_i \| \leq \gamma^\prime + hL\sum_{j< i} |a_{i,j}| \| \Psi^j - W_j \|.
\end{equation*}
We show by induction that $\| \Psi^i - W_i \| \leq \frac{\Lambda_i}{L}\gamma^\prime$: for $i=1$ one has $\| \Psi^1 - W_1 \| \leq \gamma^\prime = \frac{\Lambda_1}{L} \gamma^\prime$. For the induction step,
\begin{equation*}
    \| \Psi^i - W_i \| \leq \gamma^\prime + hL \sum_{j<i} |a_{i,j}| \frac{\Lambda_j}{L}\gamma^\prime = \frac{\gamma^\prime}{L}\big( L +hL\sum_{j<i} |a_{i,j}| \Lambda_j\big) = \frac{\Lambda_i}{L} \gamma^\prime
\end{equation*}
by definition of $\Lambda_i$. Consequently, $\| \Psi(t,Z) - W \| \leq \sum_i |b_i| \| \Psi^i - W_i \| \leq \gamma$. 

Concerning the rank bound, it suffices to show that $\rank(W_i) \leq R^i \rank(Z) + \rho_i g$ for $g = \max_j R_G(t + c_jh, \gamma^\prime)$ for $1\leq i \leq s$. We again proceed by induction: for $i=1$ one has
\begin{equation*}
    \rank(W_1) \leq R_F \rank(Z) + R_G(t, \gamma^\prime) \leq  R^1 \rank(Z) + \rho_1 g.
\end{equation*}
Using \eqref{eq: rank F bound} together with the subadditivity of the rank and setting $\tilde Y_i = Z + h \sum_{j < i}a_{i,j}W_j$ we obtain
\begin{equation*}
    \rank(W_i) \leq R_F \rank \big(\tilde Y_i \big) + R_G(t+c_ih, \gamma^\prime) \leq R_F\Big( \rank(Z) + \sum_{j < i} \rank(W_j)\Big) + g.
\end{equation*}
By the induction hypothesis, one thus has
\begin{align*}
    \rank(W_i) &\leq R_F \Big( \rank(Z) + \sum_{j<i, \ a_{i,j} \neq 0}\big( R^j \rank(Z) + \rho_j g\big)\Big) +g \\
    &\leq \rank(Z)\Big(R_F + R_F \sum_{j < i, \ a_{i,j} \neq 0} R^j \Big) + g \Big( 1 + R_F \sum_{j<i, \ a_{i,j} \neq 0} \rho_j\Big) = R^i \rank(Z) +  g \rho_i.
\end{align*}
where we used the definitions of $R^i$ and $\rho_i$ in the last step. Since $W$ approximates $\Psi(t,Z)$ with accuracy $\gamma$ and $\mathcal R_\gamma \big(\Psi(t,Z)\big)$ has minimal rank among all such approximations, it follows that
\begin{equation*}
    \rank\Big( \mathcal R_\gamma \big(\Psi(t,Z)\big) \Big) \leq \rank(W) \leq R_\Psi \rank(Z) + R^\prime \max_i R_G \Big(t+c_ih, \frac{L}{\Lambda}\gamma\Big). \qedhere
\end{equation*}
\end{proof}

\begin{remark}
The recursions stated in Proposition~\ref{prop: assumptions on F ot Psi} give the explicit bounds
\[
\Lambda\leq L\bigg(\sum_{i=1}^s |b_i|\bigg)\exp\bigg(hL \sum_{i=1}^s \sum_{j < i}|a_{i,j}|\bigg)
\]
and
    \begin{equation}\label{eq: R_psi bound}
        R_{\Psi} \leq \sum_{i=1}^s R_F\bigl(1+R_F\bigr)^{i-1} = \bigl(1+R_F\bigr)^s-1
    \end{equation}
since $R^i \leq R_F\bigl(1+R_F\bigr)^{i-1}$. In particular, $\Lambda$ can be chosen independent of $h$ since $h\leq 1$.
\end{remark}

The above bounds can be improved for specific examples.

\begin{example}
The explicit Euler scheme $\Psi(t,X)=F(t,X)$ satisfies $\Lambda=L$, $R_\Psi=R_F$ and $R^\prime = 1$, so that $R_0(t, \gamma) = R_G(t, \gamma)$. 

For the explicit midpoint method with
\[
\Psi(t,X)=F\Big(t+\frac{h}{2}, X+\frac{h}{2} F(t,X)\Big),
\]
we obtain $\Lambda= L+hL^2/2$, $R_\Psi=R_F+R_F^2$, $R^\prime = (1+ R_F)$ and
\begin{equation*}
    R_0(t,\gamma) = (1+ R_F) \max\Bigg\{ R_G\Big( t, \frac{\gamma}{1+hL/2} \Big), R_G\Big(t + \frac{h}{2}, \frac{\gamma}{1+hL/2} \Big) \Bigg\}.
\end{equation*}
Finally, for RK4, it holds that
\[
\Lambda=L+\frac{hL^2}{2}+\frac{h^2L^3}{6}+\frac{h^3L^4}{24},
\]
so that $1+h\Lambda\leq \exp(hL)$, as one would expect when looking at the ODE $\partial_tX=LX$. Moreover,
\[
R_\Psi=4R_F+3R_F^2+2R_F^3+R_F^4 \quad \text{and} \quad R^\prime = 4 + 3R_F +2R_F^2 + R_F^3,
\]
and the source term enters similarly at the times $t$, $t + h/2$ and $t+h$.

Additionally, $1+h\Lambda\leq \exp(hL)$ implies $L \gamma / \Lambda \geq e^{-hL}\gamma$, so the source is measured at essentially the accuracy $\gamma$. Note that the general value of $R_\Psi$ can be quite pessimistic: for the autonomous Sylvester right-hand side $F(X) = AX +XB$, the increment $\Psi(X)$ is a polynomial of degree $\leq s$ in the Sylvester operator applied to $X$. Thus, the column space of $\Psi(X)$ is contained in $\sum_{j=0}^sA^j \ran(X)$, so that $\rank\big(\Psi(X)\big) \leq (s+1)\rank(X)$. For RK4, one can thus take $R_\Psi = 5$ instead of $4R_F+3R_F^2+2R_F^3+R_F^4 = 52$ for $R_F = 2$.
\end{example}

In the following, we proceed by looking into the question on how one can combine explicit Runge-Kutta methods with rank reduction strategies.

\subsection{Classical step-truncation methods}\label{sec: class step-trunc}

Step-truncation methods follow a simple idea: first propagate the approximation using a standard time-stepping scheme, and then truncate the obtained approximation to a desired accuracy, in order to keep the computational costs under control. Since the prescribed parameter is the tolerance of the subsequent recompression step, such methods can be viewed as rank-adaptive procedures.

\begin{figure}[ht]
			\centering 
			\resizebox{0.5\linewidth}{!}{
				\begin{tikzpicture}[>=stealth, x=0.75cm]
					
					% axis (ending at T)
					\draw[->] (0,0) -- (10,0);
					
					% ticks
					\foreach \x/\lab in {0/{$0$},2/{$h$},4/{$2h$},%7/{$\cdots$},
                    10/{$T$}}{
						\draw (\x,0.08) -- (\x,-0.08) node[below=4pt] {\lab};
					}

					% state labels
					\node (u0)  at (0,0.33)  {$\hat X_0$};
					\node (u1)  at (2,0.33)  {$\hat{X}_1$};
					\node (u2)  at (4,0.33)  {$\hat{X}_2$};
					%\node (u3)  at (5,0.33)  {};
					\node (uN)  at (10,0.33) {$\hat{X}_N$};
					
					% propagators Psi_h (green)
					\draw[->,green!60!black,thick] (0.2,0.5)  to[bend left=35] node[above] {$\Psi$} (1.7,0.5);
					\draw[->,green!60!black,thick] (2.2,0.5)  to[bend left=35] node[above] {$\Psi$} (3.7,0.5);
					\draw[->,green!60!black,thick] (4.2,0.5)  to[bend left=35] node[above] {$\Psi$} (5.7,0.5);
					\draw[->,green!60!black,thick] (8.2,0.5)  to[bend left=35] node[above] {$\Psi$} (9.7,0.5);
					
					% self-loops (corrections)
					\draw[->, yshift=-10pt] (u0)  to[out=110,in=70,looseness=12] node[above] {$\mathcal{R}_{\delta_0}$} (u0);
                    %was [out=120,in=60,looseness=8]
					\draw[->, yshift=3pt] (u1)  to[out=110,in=70,looseness=12] node[above] {$\mathcal{R}_{\delta_1}$} (u1);
					\draw[->, yshift=3pt] (u2)  to[out=110,in=70,looseness=12] node[above] {$\mathcal{R}_{\delta_2}$} (u2);
					\draw[->, yshift=3pt] (uN)  to[out=110,in=70,looseness=12] node[above] {$\mathcal{R}_{\delta_N}$} (uN);

                    \node at (7,0.5) {$\dots$};
                    \node at (7,-0.5) {$\dots$};
				\end{tikzpicture}
			}
            \caption{Schematic representation of step-truncation methods: each iterate is produced by a Runge-Kutta step followed by a truncation to reduce the rank.} \label{fig: scheme step-trun}
		\end{figure}
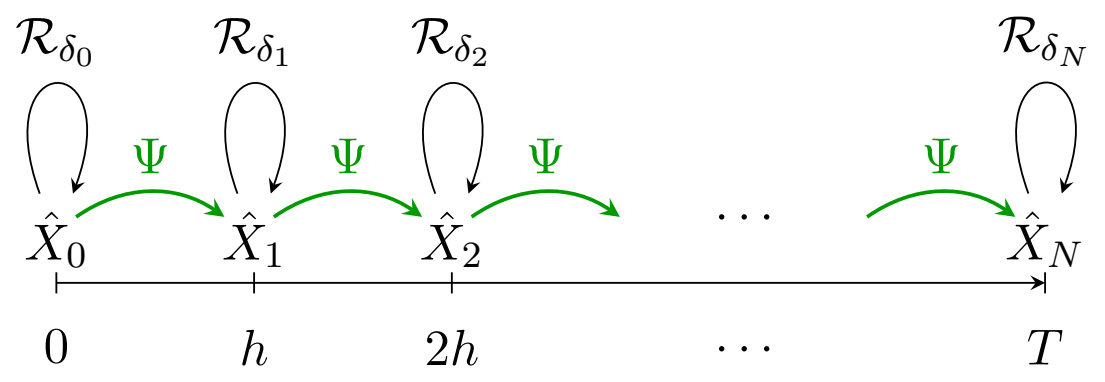

The algorithm then reads
\begin{equation}
    \hat X_{n} = \mathcal R_{\delta_{n}}\big(\hat Z_{n}\big), \quad \hat Z_{n+1} = \hat{X}_n + h \Psi(t_n, \hat{X}_n), \quad \hat Z_0 = X_0, \label{eq: step truncate} 
\end{equation}
for all $n \in \mathbb N$, which is schematically represented in Figure \ref{fig: scheme step-trun}. In the following, we give two common choices of tolerances and analyze their respective advantages and limitations.

\begin{example}[Control of the accuracy] \label{ex: classical st}
In order to preserve the convergence order of the underlying Runge-Kutta method, classical step-truncation schemes \cite{ST_2022} select the recompression tolerance at the scale of the local discretization error, that is, 
\begin{equation*}
    \delta_n = \mathcal O(h^{q+1}).
\end{equation*}
What is not clear, however, is how exactly the ranks of the computed approximations then relate to the best possible approximation ranks, since iterating the one-step bound
\begin{equation}\label{eq: rank bound u hat}
    \rank \big(\hat X_{n+1} \big) \leq \rank \big(\hat Z_{n+1} \big) \leq (1+ R_\Psi) \rank \big(\hat X_n\big) + R_0(t_n,0),
\end{equation}
which follows from \eqref{assumption-rank-Psi} with $\gamma \rightarrow 0$, yields a global rank bound that is exponential in the number of time steps. Even for homogeneous problems, where $R_0 = 0$, one only obtains
\begin{equation*}
    \rank\big(\hat X_n \big) \leq (1+R_F)^{sn}\rank \big(\hat X_0 \big).
\end{equation*}

As was already observed in \cite{BJSV26}, the analysis is considerably simplified when one compares the ranks of the computed approximation with those of the Runge-Kutta update before recompression. By definition of $\mathcal R_\delta$, the rank of $\hat X_{n+1}$ is minimal among all approximations of $\hat Z_{n+1}$ with error at most $\delta_{n+1}$. This, however, does not guarantee that the ranks are near-optimal compared to the best low-rank approximations of the exact solution. The tolerance $\delta_{n+1} = \mathcal O(h^{q+1})$ is much smaller than the distance between $\hat Z_{n+1}$ and $X(t_{n+1})$, which is typically the size of the global error $\mathcal O(h^q)$. Without further assumptions, $\hat Z_{n+1}$ and $X(t_{n+1})$ need not have a comparable decay of singular values.

To illustrate this fact, it is sufficient to consider a simple example. Let $\varepsilon \in (0,1)$, $p\in [1,\infty]$, $2 \leq k \leq d_2$, and denote by $e_i$ the canonical basis vectors. The matrix $U=e_1e_1^T$ can then be approximated by
\begin{equation*}
    W_k = e_1e_1^T + \frac{\varepsilon}{(k-1)^{1/p}} \sum_{i=2}^{k} e_ie_i^T, \qquad \| U - W_k \| = \varepsilon.
\end{equation*}
However, $U$ has only one nonzero singular value, $\sigma_1(U) = 1$, whereas the approximation $W_k$ also has the additional tail of singular values
\begin{equation*}
    \sigma_2(W_k) = \dots = \sigma_k(W_k) = \frac{\varepsilon}{(k-1)^{1/p}}.
\end{equation*}
As a consequence, for any $0 < \rho < 1$, it holds that $\rank (\mathcal R_{\rho \varepsilon}(U)) = 1$, while
\begin{equation*}
    \rank (\mathcal R_{\rho \varepsilon}(W_k)) = 1 + \lceil (1- \rho^p)(k-1) \rceil,
\end{equation*}
with the convention $\rho^p = 0$ for $p=\infty$. This grows linearly in $k$. Thus, even though $W_k$ is an accurate approximation of $U$, the two matrices do not have a comparable low-rank structure at tolerances below their distance $\varepsilon$, since even a small perturbation can create many singular values. In contrast, for $\rho \geq 1$, discarding the whole tail produces an error of $\varepsilon \leq \rho \varepsilon$, so that $\rank (\mathcal R_{\rho \varepsilon}(W_k)) \leq 1$. Thus, once the tolerance exceeds the distance between the two matrices, the ranks become comparable: this is made precise in Lemma \ref{lemma: BD recompression} below.
\end{example}

Another possibility is to choose the recompression tolerance so that one can relate the obtained ranks to the possible best approximation ranks. Therefore, we introduce the following lemma, which can be viewed as an equivalent of \cite[Lemma 20]{BDS25} for recompressions instead of soft thresholding, and is needed to control the ranks of the scheme (here $\hat X_n$) by the ranks of the recompressed solution $X(t_n)$. 
Note that this result is a special case of \cite[Lemma 2]{BachmayrDahmen15}, which studies recompression operators for tensor formats and proves near-optimality under assumptions on the rate of decrease of singular values. Indeed, it suffices to take the parameters $(\eta,\kappa_P,d,\alpha)$ of \cite[Lemma 2]{BachmayrDahmen15} equal to $(\eta,1,2, \varepsilon/\eta)$.
As the proof simplifies in the matrix case, we recall it below.

\begin{lemma}[Flexible recompression] \label{lemma: BD recompression}
Let $X, Y \in \mathbb R^{d_1 \times d_2}$. If $\|X - Y\| \leq \eta$, then
\begin{equation*}
    \rank \big(\mathcal R_{\delta} (Y)\big)\leq \rank\big(\mathcal R_{\varepsilon}(X)\big)
\end{equation*}
for every $\varepsilon>0$ and $\delta>0$ that satisfy
\begin{equation*}
    \delta \geq \eta+\varepsilon.
\end{equation*}
\end{lemma}
\begin{proof}
Observe that
\[
\|Y-\mathcal R_{\varepsilon}(X)\|\leq \|Y - X\|+\|X-\mathcal R_{\varepsilon}(X)\|\leq \eta + \varepsilon \leq \delta.
\]
Thus, $\mathcal R_{\varepsilon}(X)$ is an admissible approximation to $Y$ with tolerance $\delta$. Since $\mathcal R_{\delta}(Y)$ is chosen as an approximation of minimal rank among all approximations to $Y$ with error at most $\delta$, it follows that 
\[
\rank \big(\mathcal R_{\delta} (Y)\big)=\min_{Z:\,\| Y -Z\|\leq \delta}\rank(Z)\leq \rank \big(\mathcal R_{\varepsilon}(X)\big).\qedhere
\]
\end{proof}

\begin{example}[Control of the ranks]
    Let 
    \begin{equation*}
        \hat e_n = \|\hat X_n-X(t_n)\|
    \end{equation*}
    denote the global error of the step-truncation scheme \eqref{eq: step truncate}. By adding and subtracting $X(t_n)$ and $h \Psi(t_n, X(t_n))$, using \eqref{assumption-local-Psi} and \eqref{assumption-lipschitz-Psi}, and the fact that we employ a Runge-Kutta method of order $q$, we then know that
    \begin{align*}
        \|\hat Z_{n+1} - X(t_{n+1}) \| &= \| \hat{X}_n + h \Psi(t_n, \hat{X}_n) - X(t_{n+1})\|\\
        &\leq \| \hat X_n - X(t_n) \| + h\|\Psi(t_n, \hat X_n) - \Psi(t_n, X(t_n)) \| \\
        & \qquad \qquad \qquad \quad \, \, \,  + \| X(t_n) + h \Psi(t_n, X(t_n)) -X(t_{n+1})\| \\
        &\leq (1 + h \Lambda) \hat e_n + Ch^{q+1},
    \end{align*}
    holds by the triangle inequality.
    Setting
    \begin{equation}
        \eta_{n+1} = (1 + h \Lambda) \hat e_n + Ch^{q+1}, \label{eq: eta_n for rank control example}
    \end{equation}
    we thus obtain the upper bound
    \begin{equation}\label{eq: hat w - u}
        \|\hat Z_{n+1} - X(t_{n+1}) \| \leq \eta_{n+1}.
    \end{equation}
    Using Definition \eqref{eq: step truncate} together with the flexible recompression result of Lemma \ref{lemma: BD recompression}, for every $\varepsilon_{n+1} > 0$, if $\delta_{n+1} \geq \eta_{n+1} + \varepsilon_{n+1}$, then
    \begin{equation*}
        \rank \big(\hat X_{n+1} \big) = \rank \Big(\mathcal R_{\delta_{n+1}} \big(\hat Z_{n+1}\big)\Big) \leq \rank \Big(\mathcal R_{\varepsilon_{n+1}} \big(X(t_{n+1})\big)\Big),
    \end{equation*}
    that is, the rank of the computed approximation is controlled by the minimal rank needed to approximate the true solution $X(t_{n+1})$ with precision $\varepsilon_{n+1}$.

    The difficulty now arises when taking a closer look at the size of the total error after recompression $\hat{e}_{n+1}$, which can be bounded as
    \begin{equation*}
        \hat e_{n+1} = \|\mathcal R_{\delta_{n+1}}\big(\hat{Z}_{n+1} \big)-X(t_{n+1})\| \leq \|\mathcal R_{\delta_{n+1}}\big(\hat{Z}_{n+1} \big) - \hat{Z}_{n+1}\|
    +
    \|\hat{Z}_{n+1}-X(t_{n+1})\| \leq \delta_{n+1} + \eta_{n+1},
    \end{equation*}
    where we added and subtracted $\hat Z_{n+1}$, used the definition of $\mathcal R_{\delta_{n+1}}$, and relation \eqref{eq: hat w - u}.
    
    Then, choosing the smallest admissible recompression tolerance as
    \begin{equation*}
        \delta_{n+1} = \eta_{n+1}+\varepsilon_{n+1},
    \end{equation*}
    together with the definition of $\eta_{n+1}$ from \eqref{eq: eta_n for rank control example} yields the global error recursion
    \begin{equation}\label{eq: local err example}
        \hat e_{n+1} \leq 2(1 + h \Lambda) \hat e_n + 2Ch^{q+1} + \varepsilon_{n+1}.
    \end{equation}

    At this point, one can already recognize the problem that arises in \eqref{eq: local err example}: if we choose the recompression tolerance such that the rank control with respect to the true solution is guaranteed, then the theoretical global error bound grows exponentially in the number of time steps $n$ due to $2(1 + h \Lambda) \geq 2$, independently of the exact choice of $\varepsilon_{n+1}$.
\end{example}

As we have seen, achieving an appropriate balance between accuracy and approximation ranks is challenging for the classical step-truncation algorithm, given the existing approximation results. In the next section, we propose a strategy to overcome this issue.

\section{Step-truncation methods using additional recompressions}\label{sec: step-trunc add recom}

Motivated by the observations from \cite{BDS25}, where, in order to control the ranks of the approximation produced by an iterative scheme, one also needs to maintain rank control on the computed iterates during the algorithm, we propose to modify the algorithm by using two separate recompressions on the previously computed data. Therefore, we define the approximation scheme by
\begin{equation}\label{eq: def_scheme}
X_0 = X(t_0), \quad  X_{n+1} = \mathcal R_{\delta_n} ( X_n) + h \Psi\big(t_n, \mathcal R_{\varepsilon_n}(X_n) \big),
\end{equation}
with recompression tolerances $\delta_n, \varepsilon_n \geq 0$ to be specified in what follows for all $0\leq n < N$.

\subsection{Choices of tolerances, periodic rank resets, and global accuracy result} \label{sec: global accuracy 1}

We start by specifying the exact tolerances. As was shown in \cite{BJSV26}, at least for the norm preserving Schr\"odinger equation, it is possible to obtain quasi-optimal rank bounds, which depend on the time step size $h$ and final time $T$ proportionally to $T/h$. In the following, we introduce a possible strategy in order to also avoid the linear dependence of the ranks on the final time $T$.

Let $\xi_n$, $0 \leq n \leq N$, be defined as in \eqref{eq: thm 1 opt ranks}. We show below that $\| X_n - X(t_n) \| \leq \xi_n$. Let further $M = 1/h \in \mathbb N$. The recompression tolerances are then defined as
\begin{equation}
    \varepsilon_n = 2 \xi_n \quad \text{and} \quad \delta_n = \begin{cases}
    	    \xi_0,\quad & \text{if } n=0,\\
            2 \xi_n, \quad &\text{if } n \in M \mathbb N = \{M, 2M, \dots \},\\
            2h \xi_n, \quad & \text{otherwise}.
        \end{cases} \label{eq: tolerances adaptive}
\end{equation}

This particular choice of the first recompression tolerance $\delta_n$ allows us to reset the approximation ranks accumulated over the last $1/h$ steps back to optimal ones, see the rank bound stated in Theorem \ref{thm:main}. For the error propagation, this choice does not pose a problem since this reset happens only every $1/h$ time steps and thus one is able to avoid an exponential blow-up of the error in $n$. Overall, we thus obtain the modified step-truncation algorithm presented in Algorithm \ref{alg:step_trunc_reset}. Note that the algorithm is formulated using variable tolerances $\xi_1, \dots, \xi_{N-1}$, as commonly used in numerical implementations.

\begin{algorithm}[ht]
\caption{Step truncation with separate recompressions and periodic rank resets}
\label{alg:step_trunc_reset}
\begin{algorithmic}[1]
\Require step size $h\in(0,1]$ with $M=1/h\in\mathbb N$, number of steps $N$, initial data $X_0=X(0)$, constants $C$ and $\Lambda$, tolerances $\xi_0,\dots,\xi_N$
\For{$n=0,\dots,N-1$}
    \State $t_n \gets nh$
    \If{$n=0$}
        \State $\delta_n \gets \xi_0$
    \ElsIf{$n \bmod M = 0$} \Comment{periodic rank reset}
        \State $\delta_n \gets 2\xi_n$
    \Else
        \State $\delta_n \gets 2h\xi_n$
    \EndIf
    \State $X_{n+1} \gets \mathcal R_{\delta_n}(X_n) + h\,\Psi\bigl(t_n,\mathcal R_{2\xi_n}(X_n)\bigr)$
\EndFor
\State \Return $X_0,\dots,X_N$
\end{algorithmic}
\end{algorithm}

We start the proof of the main theorem by first addressing the stated accuracy result.

\begin{proof}[Proof of the accuracy result in Theorem \ref{thm:main}]
 
We proceed by induction on $n$. For the initial datum it holds that $e_0 =0\leq \xi_0$.

Now let $n\geq 0$ such that $e_n \leq \xi_n$. By assumption~\eqref{assumption-local-Psi},
\[
\|X(t_{n+1})-X(t_n)-h\Psi(t_n,X(t_n))\|\leq C h^{q+1}\leq h\xi_n.
\]
Using triangle inequalities, we obtain
\[
e_{n+1} =  \| X_{n+1} - X(t_{n+1})\| \leq Ch^{q+1}+e_n+\delta_n+h\Lambda(e_n+\eps_n).
\]
If $n=0$, as $e_0=0$, we have
\[
e_1\leq Ch^{q+1}+\xi_0+h\Lambda \eps_0=(1+h+2h\Lambda)\xi_0\leq \xi_1.
\]
If $n\in M\N$, observing that $Ch^{q+1}\leq h\xi_n$, we get
\[
e_{n+1}\leq (3+h+3h\Lambda)\xi_n\leq 3\exp(3(\Lambda+1)h)\xi_n=\xi_{n+1}.
\]
Otherwise, as $\delta_n=2h\xi_n$, we obtain
\[
e_{n+1}\leq (1+3h+3h\Lambda)\xi_n\leq \exp(3(\Lambda+1)h)\xi_n \leq \xi_{n+1}. \qedhere
\]
\end{proof}

\begin{remark}[Optimality of the error bound] \label{rem: exp T dependence}
    The exponential dependence on the final time~$T$ in Theorem~\ref{thm:main} is, in fact, not an artifact of the proof, but generally unavoidable for non-contractive evolution equations. Indeed, if the flow is only Lipschitz continuous, then perturbations are propagated according to a Gr\"onwall-type estimate. More precisely, for two solutions $X(t)$ and $Y(t)$ of \eqref{eq: ODE} with different initial data $X(0)$ and $Y(0)$, respectively, it holds that
    \begin{equation*}
        \|X(t)-Y(t)\|
        \leq
        e^{L t}\|X(0)-Y(0)\|.
    \end{equation*}
    
    Thus, an initial perturbation may be amplified by a factor of order $e^{L T}$. Furthermore, local perturbations caused by the 
    time-discretization and the recompression might cause an error of a similar factor. Thus, no error bound which is uniform for all Lipschitz, non-contractive right-hand sides can avoid an exponential dependence on $T$. 
    Moreover, a classical error analysis of the scheme obtained by simply iterating time steps with increment $\Psi$ would give a bound of the form $Ch^q\exp((\Lambda+1)t_n)$. In Theorem~\ref{thm:main}, the error bound is similar, with the constant $\Lambda+1$ replaced by $3\Lambda+3+\ln(3)$.
    
    A better dependence on $T$ is only possible under additional stability
    assumptions, for example if the flow is contractive or dissipative. Without such additional structures, the exponential-in-$T$ growth is sharp.
\end{remark}

\subsection{Obtaining quasi-optimal rank bounds}

Next, we prove that by choosing the tolerances according to \eqref{eq: tolerances adaptive} we also maintain rank control on the computed approximations. Therefore, we follow the general strategy developed in \cite{BDS25}. Since we do not apply any soft thresholding operator in the algorithm, we can further simplify the proof of rank bounds.

We start by providing an equivalent of \cite[Lemma 23]{BDS25} for an arbitrary Schatten $p$-norm, allowing us to lower the truncation threshold by a factor $\theta$, up to an increase in the ranks.
Recall that the recompression operator $\mathcal R$ is defined in \eqref{eq: R_delta} with respect to the norm $\|\cdot\|=\|\cdot\|_p$ from \eqref{eq:schatten}. 
\begin{lemma} \label{lemma: rank h xi}
For any $\xi>0$, $\theta \in (0,1)$, $p \in [1,\infty]$ and $X\in \R^{d_1\times d_2}$ of rank $R=\rank(X)$, it holds that
\[
\rank \big(\mathcal R_{\theta\xi}(X) \big)\leq \lceil \theta^p r+(1-\theta^p)R\rceil,
\]
where $r=\rank \big(\mathcal R_{\xi} (X) \big)$. Note that we use the convention $\theta^\infty = 0$ for $p = \infty$.
\end{lemma}

\begin{proof}
If $r = R$, the statement of the lemma follows immediately. Therefore, assume that $r < R$. Denote the singular values of $X$ sorted in decreasing order by $\sigma_1\geq \dots \geq \sigma_R$, then the definition of $\mathcal{R}_{\xi}$ and the choice of $r$ imply that
\[
\sum_{k=r+1}^R \sigma_k^p\leq \xi^p.
\]
Setting $r'=\lceil \theta^p r+(1-\theta^p)R\rceil\geq r$, we again notice that the case $r' = R$ directly yields the statement of the lemma and we can further assume $r' < R$. As the $\sigma_k$ are non-increasing, the average of the last $R-r'$ singular values to the power $p$ is not larger than the average of the last $R-r$, that is,
\[
\frac{1}{R- r'}\sum_{k=r'+1}^R \sigma_k^p\leq \frac{1}{R-r}\sum_{k=r+1}^R \sigma_k^p.
\]
Combining these two inequalities yields
\[
\sum_{k=r'+1}^R \sigma_k^p \leq \frac{R-r'}{R-r}\xi^p\leq  \frac{R-(\theta^p r+(1-\theta^p)R)}{R-r}\xi^p=\theta^p\xi^p.
\]
These computations show that $\rank \big(\mathcal R_{\theta\xi}(X)\big)\leq r'$, which proves the claim.
\end{proof}
Note that the summability parameter $p$ has an influence on the rank bound of Lemma~\ref{lemma: rank h xi}. In particular, the result is stronger for $p=1$, corresponding to the nuclear norm $\|\cdot\|_1$. However, the choice of $p$ also depends on the application: when discretizing the Schr\"odinger equation, the natural choice is the Frobenius norm $\|\cdot\|_2$, which is equivalent to the $L^2$ norm of the wave function \cite{BDS25, BJSV26}.

Motivated by \cite[Proposition 25]{BDS25}, Lemmas~\ref{lemma: BD recompression} and \ref{lemma: rank h xi} can now be used in the following result in order to establish rank control over multiple intervals in the form of a recursive rank relation.

\begin{lemma} \label{lem: recursion}
    Let $\rho_n = \rank\big( \mathcal R_{\delta_n}(X_n)\big)$ denote the rank of the recompressed numerical approximation $X_n$. Then $\rho_0 = \bar r_0$ and $\rho_n \leq \bar r_n$ for $n \in M \mathbb N$. For $0 \leq n \leq N-1$ with $n+1 \notin M \mathbb N$ it holds that
    \begin{equation*}
        \rho_{n+1} \leq 1 + \bar r_{n+1} + e^{-(h/3)^p}\big( \rho_n + R_\Psi \bar r_n + R_0(t_n, \xi_n) \big).
    \end{equation*}
    Moreover, for $0 \leq n \leq N-1$,
    \begin{equation*}
        \rank(X_{n+1}) \leq \rho_n + R_\Psi \bar r_n + R_0(t_n,0).
    \end{equation*}
    Note that $R_{\Psi}$ and $\bar r_n$ are defined in \eqref{assumption-rank-Psi} and \eqref{eq: thm 1 opt ranks}, respectively.
\end{lemma}

\begin{proof}

Recall that the optimal approximation ranks satisfy $\bar r_n = \rank \Big( \mathcal R_{\xi_n}\big(X(t_n)\big) \Big)$. For the initial data we have $\delta_0 = \xi_0$ and $X_0 = X(0)$, so $\rho_0 = \rank \big( \mathcal R_{\xi_0}(X(0)) \big) = \bar r_0$.

Since $e_n = \| X_{n} - X(t_{n}) \| \leq \xi_{n}$, Lemma \ref{lemma: BD recompression} with $\eta = \varepsilon = \xi_n$ yields $\rank \big(\mathcal R_{2 \xi_{n}}(X_{n})\big) \leq \bar{r}_{n}$ for all $0 \leq n \leq N$. At reset steps, using $\delta_n = 2 \xi_n$, we then obtain $\rho_n \leq \bar r_n$.

By subadditivity and \eqref{assumption-rank-Psi} with $\gamma \rightarrow 0$, it then holds that
\begin{equation*}
    \rank(X_{n+1}) \leq \rank\big( \mathcal R_{\delta_n}(X_n)\big) +\rank \Big( \Psi \big(t_n, \mathcal R_{2 \xi_n}(X_n) \big) \Big)\leq \rho_n + R_\Psi \bar r_n + R_0(t_n,0),
\end{equation*}
which corresponds to the last bound stated in the lemma.

Lastly, for non-reset steps, we first define the approximations
\begin{equation*}
    u_n = \mathcal R_{\delta_n}(X_n) + h W_n, \qquad W_n = \mathcal R_{\xi_n} \Big( \Psi \big(t_n, \mathcal R_{2 \xi_n}(X_n) \big) \Big),
\end{equation*}
for $0 \leq n \leq N-1$, satisfying $\| X_{n+1} - u_n \| \leq h \xi_n$. Using \eqref{assumption-rank-Psi}, we obtain $\rank(W_n) \leq R_\Psi \bar r_n + R_0(t_n, \xi_n)$ as well as 
\begin{equation}
    \rank(u_n) \leq \rho_n + R_\Psi \bar r_n + R_0(t_n, \xi_n). \label{eq: rank rec 3}
\end{equation}
Since $\xi_n \leq \xi_{n+1}$, we can use Lemma \ref{lemma: BD recompression} with the tolerances $\eta = h \xi_n$, $\varepsilon = h \xi_{n+1}$ and $\delta = 2h \xi_{n+1}$ to arrive at
\begin{equation}
    \rho_{n+1} = \rank\big( \mathcal R_{2 h \xi_{n+1}}(X_{n+1})\big) \leq \rank \big( \mathcal R_{h \xi_{n+1}}(u_n) \big). \label{eq: rank rec 1}
\end{equation}
Moreover, since also $\|u_n - X(t_{n+1}) \| \leq h \xi_n + \xi_{n+1} \leq 2 \xi_{n+1}$, we can again use Lemma \ref{lemma: BD recompression} with $\eta = 2 \xi_{n+1}$, $\varepsilon = \xi_{n+1}$ and $\delta = 3 \xi_{n+1}$ to obtain
\begin{equation}
    \rank\big( \mathcal R_{3 \xi_{n+1}}(u_n) \big) \leq \rank \Big( \mathcal R_{\xi_{n+1}}\big( X(t_{n+1})\big) \Big) = \bar r_{n+1}. \label{eq: rank rec 2}
\end{equation}
The full chain for a non-reset step then reads
\begin{align*}
    \rho_{n+1} &\leq \rank \big( \mathcal R_{h \xi_{n+1}}(u_n) \big) = \rank \big( \mathcal R_{3 \theta \xi_{n+1}}(u_n) \big) \\
    &\leq \left\lceil \theta^p \rank\big( \mathcal R_{3 \xi_{n+1}}(u_n) \big) + (1- \theta^p) \rank(u_n) \right\rceil \\
    &\leq \left\lceil \theta^p \bar r_{n+1} + (1- \theta^p)(\rho_n + R_\Psi \bar r_n + R_0(t_n, \xi_n)) \right\rceil,
\end{align*}
for $\theta = h/3$, where in the first step we used \eqref{eq: rank rec 1}, in the second inequality Lemma \ref{lemma: rank h xi} and lastly \eqref{eq: rank rec 2} together with \eqref{eq: rank rec 3}. Finally, using $\theta^p \leq 1$, $1-\theta^p \leq e^{-\theta^p}$, $\lceil x\rceil \leq x +1$ and inserting the corresponding value of $\theta$ yields the claimed recursive bound for non-reset steps.
\end{proof}

We are now in position to also prove the global rank bound stated in Theorem \ref{thm:main}.

\begin{proof}[Proof of rank bounds in Theorem \ref{thm:main}]
For $n\in \N$, let $n'=\lfloor t_n\rfloor/h$ be the largest multiple of $M$ smaller or equal to $n$ and set $a = e^{-(h/3)^p}$. By induction on Lemma \ref{lem: recursion} and $\rho_{n^\prime} \leq \bar r_{n ^\prime}$, it holds that
\[
\rho_{n}\leq \sum_{k=n'}^n a^{n-k} \bar r_k + \sum_{k = n^\prime}^{n-1} a^{n-k} \big( R_\Psi \bar r_k + R_0(t_k, \xi_k) \big) + (n - n^\prime).
\]
Inserting this into the corresponding rank bound for $X_{n+1}$ then yields
\begin{equation*}
    \rank(X_{n+1}) \leq \sum_{k = n^\prime}^n a^{n-k} \big((1 + R_\Psi) \bar r_k + R_0(t_k, \xi_k) \big) + (n - n^\prime) + R_0(t_n,0).
\end{equation*}
Taking the maximum over the summands and using that $n - n^\prime +1 \leq 1/h$ and $h^p(n-k) = h^{p-1}(t_n - t_k)$, we obtain the stated uniform rank bound.

The sharper bound for reset times follows directly from combining the second bound of Lemma \ref{lem: recursion} with $\rho_n \leq \bar r_n$ for reset times $t_n$.
\end{proof}

In summary, we showed that for our proposed scheme, the ranks of the obtained approximations do not grow exponentially in the total number of time steps: at each reset step they can be bounded by the corresponding optimal approximation ranks (with respect to the current error), and between those resets they grow by at most a factor of order $1/h$, which we will address in the following section. Note that the source contributes only $R_0(t_k, \xi_k)$ per step, that is, the rank required for an approximation to the current error, and the exact rank $R_0(t_n,0)$ appears only once.

\begin{remark}[Alternative placement of thresholds] \label{rem: alt thresh}

The proposed scheme can also be written in a form closer to classical step truncation, with one recompression after each time step and just an additional inner recompression:
\begin{equation}\label{eq: def_scheme 2}
        \tilde Y_0 = \mathcal R_{\delta_0}\big(X(0)\big), \qquad  Y_{n+1} =  \tilde Y_n + h \Psi \big(t_n, \mathcal R_{\varepsilon_n}( \tilde Y_n) \big), \qquad \tilde Y_{n+1} = \mathcal R_{\delta_{n+1}}(Y_{n+1}),
\end{equation}
with $\varepsilon_n = 4 \xi_n$ and $\delta_n$ as in \eqref{eq: tolerances adaptive}. The inner tolerance is doubled compared to \eqref{eq: tolerances adaptive}, since $\tilde Y_n$ is generally only within $3 \xi_n$ of $X(t_n)$.

In fact, this is essentially the same algorithm as \eqref{eq: def_scheme}: substituting $\tilde Y_n = \mathcal R_{\delta_n}(Y_n)$ yields
\begin{equation*}
        Y_{n+1} = \mathcal R_{\delta_n}(Y_n) + h \Psi \Big(t_n, \mathcal R_{4 \xi_n}\big(\mathcal R_{\delta_n}(Y_n)\big) \Big), \qquad Y_0 = X(0).
\end{equation*}
The two methods differ only in two places: the propagated argument of $\Psi$ and which iterate is reported. Since $\mathcal R_{4 \xi_n}\big(\mathcal R_{\delta_n}(Y_n)\big)$, as a composition, is a truncation of $Y_n$ at a tolerance between $4\xi_n$ and $4 \xi_n + \delta_n$, the variant uses a coarser inner tolerance. Additionally, the main scheme returns $X_n$ before the recompression $\mathcal R_{\delta_n}$, while the alternative variant returns the recompressed $\tilde Y_n$.

As a consequence, Theorem \ref{thm:main} carries over: the proofs use the second argument $T_n$ of $\Psi$ only through $\rank(T_n) \leq \bar r_n$ and $\|T_n - X(t_n) \| \leq 3 \xi_n$. For the variant, Lemma \ref{lemma: BD recompression} still yields $\rank \Big( \mathcal R_{4 \xi_n}\big(\mathcal R_{\delta_n}(Y_n)\big)\Big) \leq \bar r_n$ and $\| \mathcal R_{4 \xi_n}\big(\mathcal R_{\delta_n}(Y_n)\big) - X(t_n) \| \leq 7 \xi_n$, hence the error bound $\| Y_n - X(t_n)\| \leq \xi_n$ and the global rank bound for the update $Y_{n+1}$ hold with $\xi_n = Ch^q 3^{j_n}e^{(3+7 \Lambda)t_n}$. In addition, it also holds that
\begin{equation*}
    \big\| \tilde Y_n - X(t_n) \big\| \leq 3 \xi_n, \qquad \rank \big( \tilde Y_{n+1} \big) \leq \rank(Y_{n+1}).
\end{equation*}
The coarser inner tolerance thus leads (in theory) to a faster growth of the error bound. We compare both schemes numerically in Section \ref{sec: num exp}.
\end{remark}

\begin{remark}[Comparison to the rank bounds of {{\cite{BDS25}}} and {{\cite{BJSV26}}}]

In the matrix case, Theorem \ref{thm:main} improves on the rank bounds available for iterative thresholding schemes in several aspects.

Compared with \cite[Theorem 27]{BDS25}, whose rank bounds scale as $h^{-3}$, the ranks in Theorem \ref{thm:main} are near-optimal up to a factor $h^{-1}$. This is due to the explicit time stepping, which does not require a fixed-point iteration for an implicit system.

The bound for the exact solution in \cite[Corollary 5.8]{BJSV26} also contains one factor $h^{-1}$, but it is multiplied by the final time $T$, since ranks accumulate over all time steps. Moreover, quasi-optimality is established for algebraic and exponentially decaying singular values. Our present work compares the computed ranks directly with the best approximation ranks of the exact solution at the target accuracy, without any assumptions on the decay of singular values. Additionally, using the periodic rank resets, our rank bound is independent of the final time $T$; this improvement could also be obtained for the implicit schemes of \cite{BDS25} and \cite{BJSV26}, by applying a similar periodic rank reset. In all cases, the error bounds grow exponentially in $T$ when ranks are compared with the exact solution. Error bounds that are linear in $T$ are obtained for the linear Schr\"odinger equation in \cite[Theorem~4.6]{BJSV26}, but with ranks compared to those of the perturbed discrete fixed points rather than the exact solution.

Moreover, \cite{BDS25} and \cite{BJSV26} approximate linear Schr\"odinger equations using implicit collocation methods combined with iterative soft thresholding of the iterates, whereas Theorem \ref{thm:main} applies to explicit Runge-Kutta methods for general right-hand sides satisfying assumptions \eqref{eq: F smooth} -- \eqref{eq: rank F bound}, including source terms. The use of implicit methods was motivated by the possibility of integrating over longer time intervals, and refining the scheme by raising the order, making these approaches closer in spirit to space-time formulations. Our strategy lies on the other end of the spectrum, keeping the order fixed and refining by decrease of $h$.
Finally, we argue that the remaining factor $h^{-1}$ cannot be avoided in general, as illustrated in the following section.
\end{remark}

\section{Sharpness of the obtained rank bounds}\label{sec: optimality of the bounds}

In this section, we discuss the dependence of the rank bound in Theorem \ref{thm:main} on the time step $h$: in Section \ref{sec: canc smal}, we show that the factor $1/h$ is generally sharp and cannot be improved for Algorithm \ref{alg:step_trunc_reset}. Sections \ref{sec: wrong sing} and \ref{sec: miss canc} illustrate two further difficulties that prevent $h$-uniform rank bounds for related strategies, namely recompression tolerances at the scale of the local error, as in classical step truncation methods, and error analyses based on the accumulation of local errors.

In the following examples, we take $X(0) = 0$ and consider diagonal right-hand sides $F(t)$ that are Lipschitz continuous with respect to time with constant $L_t$. We restrict ourselves to the unit time interval $[0,1]$ with $M = N = 1/h$ and denote by $E_k$ the diagonal matrix with a single entry 1 at position $(k,k)$; the dimensions $d_1 \geq d_2$ are then assumed to be large enough to accommodate all entries.

For the explicit Euler scheme, we have $\Psi(t,X) = F(t)$, such that \eqref{assumption-lipschitz-Psi} holds for any $\Lambda > 0$, \eqref{assumption-rank-Psi} with $R_\Psi = 0$ and $R_0(t, \gamma) = \rank \Big( \mathcal R_{\gamma}\big( F(t) \big) \Big)$ and \eqref{assumption-local-Psi} is satisfied for $q=1$ and $C=L_t/2$. Note that it is also possible to construct such examples with more regular $F$, at the price of additional technical difficulties.

In all three examples, the exact solution satisfies $\| X(t_n)\|_F \leq h$ at all grid points for $h$ sufficiently small, which is why all best approximation ranks vanish. In particular, the trivial approximation $X_n = 0$ has then rank zero and error at most $h$. However, this approximation is not produced by a numerical time-stepping method: it ignores the right-hand side $F$ and fully relies on knowing the solution in advance.

\subsection{The recompression could act on the wrong singular values} \label{sec: wrong sing}

We start by considering the classical step-truncation method as in Example \ref{ex: classical st} with the explicit Euler scheme and the tolerance $\delta_n = h^2$ at the scale of the local error, that is,
\begin{equation*}
    X_{n+1}=\mathcal R_{h^2} \big( X_n+hF(t_n) \big), \qquad X_0 = 0.
\end{equation*}

The right-hand side is chosen such that
\begin{equation*}
    F(t_n) = \frac h2 E_1 + \frac h2 E_2 + hE_{n+3}, \qquad 0 \leq n \leq N,
\end{equation*}
and extended to $[0,1]$ by linear interpolation in time: using the nodal hat functions
\begin{equation*}
    \rho_m(t) = \max \Big\{ 1- \frac{|t - t_m|}{h}, 0 \Big\}, \qquad 1 \leq m \leq N,
\end{equation*}
and the boundary-corrected hat
\begin{equation*}
    \rho_0(t) = \begin{cases}
        1 + t/h, \quad &0 \leq t \leq h/2, \\
        3(1- t/h), \quad &h/2 \leq t \leq h, \\
        0, \quad &\text{otherwise},
    \end{cases}
\end{equation*}
which satisfy $\rho_0(0) = 1$ and $\int_0^h \rho_0(s) \, ds = h$, we set
\begin{equation*}
    F(t) = \frac{h}{2}E_1 + \frac{h}{2}E_2 + h \sum_{m=0}^{N} \rho_m(t)E_{m+3}.
\end{equation*}
Since at any time at most two hat functions are nonzero, we have $\rank \big( F(t) \big) \leq 4$ and $F$ is Lipschitz continuous in time.

The exact solution at the grid points is then given by
\begin{equation*}
    X(t_n)= \frac{ht_n}{2} (E_1 + E_2) +  h^2\sum_{m=0}^{n-1} E_{m+3} + \frac{h^2}{2}E_{n+3}, \qquad 1 \leq n \leq N.
\end{equation*}
Using $t_n \leq 1$ and $n \leq 1/h$, it follows that
\begin{equation*}
    \| X(t_n) \|_F^2 = \frac{h^2t_n^2}{2} + h^4 \Big( n + \frac 14 \Big) \leq \frac{h^2}{2} + h^3 + \frac{h^4}{4} \leq h^2, \qquad \text{for} \quad  h \leq \frac 14,
\end{equation*}
and thus $\mathcal R_h \big( X(t_n) \big) = 0$. At the global accuracy $h$, the exact solution therefore has best approximation rank zero.

\begin{figure}[ht]
    \centering
    \begin{tikzpicture}[x=0.8cm,y=0.55cm,>=Latex]

% Axes
\draw[->] (0,0) -- (9.2,0) node[right] {time step $n$};
\draw[->] (0,0) -- (0,8.4) node[above] {diagonal index};

% Time grid
\foreach \n in {0,...,8} {
    \draw[gray!30] (\n,0) -- (\n,8);
    \node[below] at (\n,0) {\scriptsize $\n$};
}

% Diagonal index labels
\foreach \j in {1,...,8} {
    \draw[gray!20] (0,\j) -- (9,\j);
    \node[left] at (0,\j) {\scriptsize $\j$};
}

% Red crosses: fixed small increments h^2/2 removed at every step
\foreach \n in {0,...,7} {
    \foreach \j in {1,2} {
        \draw[red!75,line width=0.8pt]
            (\n+1,\j) +(-0.12,-0.12) -- +(0.12,0.12);
        \draw[red!75,line width=0.8pt]
            (\n+1,\j) +(-0.12,0.12) -- +(0.12,-0.12);
    }
}

% Blue persistent h^2 modes
\foreach \j in {3,...,8} {
    \pgfmathtruncatemacro{\birth}{\j-2}
    \draw[blue!70,line width=1.7pt] (\birth,\j) -- (8.5,\j);
    \fill[blue!70] (\birth,\j) circle (2pt);
}

% Labels
\node[red!75,align=center] at (10.2,1.5)
    {\scriptsize $\frac{h^2}{2}$-increments\\[-1mm]
     \scriptsize are removed};

\node[blue!70,align=center] at (10.6,7.8)
    {\scriptsize $h^2$-modes survive\\[-1mm]
     \scriptsize and accumulate};

% Brace indicating rank growth
\draw[decorate,decoration={brace,mirror,amplitude=5pt}]
    (8.7,3) -- (8.7,8)
    node[midway,right=6pt,align=left]
    {\scriptsize $\rank(X_n)\sim n = \mathcal O(h^{-1})$\\[-1mm]
     \scriptsize at $T=1$};

\end{tikzpicture}
\caption{Schematic representation of the example in Section \ref{sec: wrong sing}. Each time step creates a new $h^2$-mode which is kept by the recompression $\mathcal R_{h^2}$, while the two smaller $\frac{h^2}{2}$-increments are removed. Thus the numerical rank grows like $n$, although the exact solution has best approximation rank zero at the global tolerance $h$.}
\label{fig: counter ex 1}
\end{figure}
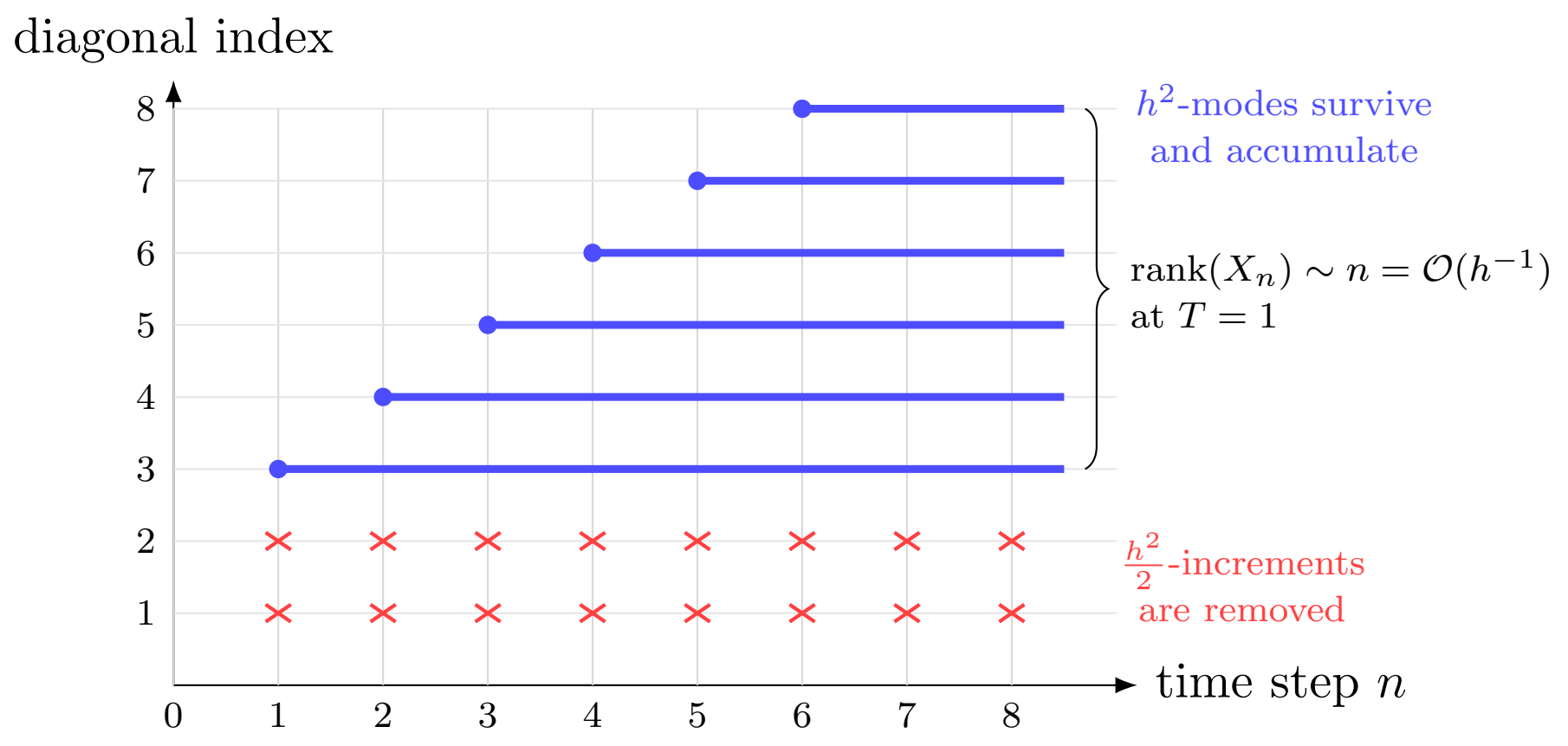

At the same time, the numerical solution satisfies by induction
\begin{equation*}
    X_n =h^2 \sum_{m = 0}^{n-1} E_{m+3}, \qquad 0 \leq n \leq N.
\end{equation*}
Indeed, if $X_n$ has this form, then $X_n + hF(t_n)$ is diagonal with $n+1$ entries of size $h^2$ at positions $3,4, \dots, n+3$ and two entries of size $h^2/2$ at the first two positions. The subsequent recompression by $h^2$ will always cancel the first two values $\frac{h^2}{2}$, since discarding in addition any of the $h^2$ entries would already create an error larger than $h^2$. Consequently, it holds that
\begin{equation*}
    \rank(X_n) = n \quad \text{and thus also} \quad \rank(X_N) = \frac 1h,
\end{equation*}
with $\| X_N - X(1) \|_F = (h^2/4 + h^2/4 + h^4/4)^{1/2} \leq h$, see also Figure \ref{fig: counter ex 1} for an illustration of the argument. 

Note that this mechanism does not affect Algorithm \ref{alg:step_trunc_reset}: at the grid points $t_n$, it holds that 
\begin{equation*}
    F(s) - F(t_n) = h(\rho_n(s)-1)E_{n+3} + h \rho_{n+1}(s)E_{n+4}, \qquad s \in [t_n , t_{n+1}],
\end{equation*}
since at most two hat functions are nonzero in $[t_n , t_{n+1}]$. The local error is thus given by
\begin{equation*}
    X(t_n+h) - X(t_n) - hF(t_n) = \int_{t_n}^{t_n+h} F(s) - F(t_n) \, ds = \frac{h^2}{2}(E_{n+4} - E_{n+3}), 
\end{equation*}
with norm $\big\|\frac{h^2}{2}(E_{n+4}-E_{n+3}) \big\|_F = h^2/\sqrt{2}$, which implies that Assumption \eqref{assumption-local-Psi} can only hold with constant $C \geq 1/\sqrt{2}$. The fine recompression tolerance can thus be lower bounded through $\delta_n = 2h \xi_n \geq 2Ch^2 \geq \sqrt{2} h^2$. The entries created in one time step have total norm $\|hF(t_n)\|_F = \sqrt{3/2}h^2 \leq \sqrt{2}h^2$, so each increment is discarded entirely in the next time step. The scheme thus satisfies $X_n = hF(t_{n-1})$ with $\rank(X_n) = 3$ over the entire time interval.

One possible way to address this particular mechanism is to truncate the forcing contribution $F(t_n)$ itself, or to recompress the component of $X_n$ orthogonal to the range of $F(t_n)$. In this way, one tries to avoid accumulating rank in directions that are introduced only by the current forcing term. However, such local modifications are not sufficient in general, as the following example shows.

\subsection{The recompression may fail to remove small singular values} \label{sec: canc smal}

The following construction shows that the factor $1/h$ in Theorem \ref{thm:main} cannot be improved for Algorithm~\ref{alg:step_trunc_reset}. Here, many pulses in different directions create singular values that are individually too large to be removed by the fine recompression $\mathcal R_{2h \xi_n}$, although their total contribution is negligible at the target accuracy.

Let $H(t) = \min\{t,1-t \}_+$ denote the hat function supported in $[0,1]$, which satisfies $\| H\|_\infty = \frac 12$, $\int_0^1 H(s) \ ds = \frac 14$ and has Lipschitz constant 1. For parameters $\alpha, \beta \in (0,1)$ with $\beta < \alpha$, we set $J_h = \lfloor h^{- \beta} \rfloor$ and 
\begin{equation}
    F(t) = \sum_{j=0}^{J_h -1} h^\alpha H(h^{- \alpha}t-j) E_{j+1}. \label{eq: ex 2 F}
\end{equation}
The $(j+1)$-th diagonal entry of $F$ thus has support $[h^\alpha j, h^\alpha(j+1)]$, so that $\rank\big(F(t)\big) \leq 1$ at all times and all pulses are contained in the interval $[0, h^{\alpha - \beta}] \subset [0,1)$. By construction, the derivative $F^\prime(t)$ has at most one nonzero entry equaling $\pm 1$, so that $F$ is also Lipschitz continuous in time with constant $L_t = 1$ in every Schatten $p$-norm. Since moreover
\begin{equation*}
     \Bigg \| \int_t^{t+h} F(s) - F(t) \, ds \Bigg\|_p \leq \int_t^{t+h} (s-t) \, ds = \frac{h^2}{2},
\end{equation*}
Assumption \eqref{assumption-local-Psi} is satisfied for $C = \frac 12$ and $q=1$ for the explicit Euler method. 

Integrating over time, the exact solution $X(t)$ is diagonal with entries being bounded through
\begin{equation*}
    0 \leq X_{j+1}(t) \leq h^\alpha \int_{h^\alpha j}^{h^\alpha(j+1)} H(h^{-\alpha}s-j) \, ds = \frac{h^{2\alpha}}{4},
\end{equation*}
and since at most $J_h$ are nonzero, we obtain
\begin{equation}
    \| X(t) \|_p \leq J_h^{1/p}\frac{h^{2\alpha}}{4} \leq \frac{h^{2 \alpha - \beta/p}}{4}, \qquad t \in [0,1]. \label{eq: ex 2 bound}
\end{equation}

\begin{prop}[Sharpness of the factor $1/h$] \label{prop: h sharpness}
    Let $p \in [1, \infty]$, $\frac 12 < \alpha < 1$ and $0 < \beta < \min\{ \alpha, p(2\alpha - 1) \}$, with the convention $p(2 \alpha -1) = \infty$ for $p = \infty$. Consider Algorithm \ref{alg:step_trunc_reset} with the explicit Euler method for $F$ as in \eqref{eq: ex 2 F}, $C = \frac 12$ and some $\Lambda > 0$. Then, for all sufficiently small $h$, it holds that
    \begin{enumerate}
        \item[(i)] $\bar r_n = 0$ for $0 \leq n \leq M$, so that Theorem \ref{thm:main} yields $\rank(X_{n+1}) \leq \frac2h +1$, $0 \leq n < M$.
        \item[(ii)] $\rank(X_n) \geq J_h = \lfloor h^{-\beta} \rfloor$ for all $n$ satisfying $J_hh^\alpha \leq t_n \leq 1$.
    \end{enumerate}
    In particular, the factor $1/h$ in Theorem \ref{thm:main} cannot be replaced by $h^{-\beta^\prime}$ for any $\beta ^\prime <1$.
\end{prop}

\begin{figure}[ht]
    \centering
    \begin{tikzpicture}[x=1.05cm,y=1cm,>=Latex]

    % Top plot: non-overlapping hats
    \draw[->] (0,0) -- (7.2,0) node[right] {time $t$};
    \draw[->] (0,0) -- (0,1.25) node[above] {$F_j(t)$};
    
    \foreach \j/\col in {0/blue,1/teal,2/violet,3/orange,4/green!60!black,5/red} {
        \draw[\col,line width=1.2pt]
            (\j,0) -- (\j+0.5,0.9) -- (\j+1,0);
        \node[below] at (\j+0.5,0) {\scriptsize $j = \j$ };
    }
    
    \node[align=center] at (3.5,1.35)
        {\small non-overlapping pulses: $\rank(F(t))\leq 1$};
    
    % Dashed vertical line showing only one active pulse
    \draw[dashed,gray] (5.8,0) -- (5.8,1.0);
    \node[gray,align=center,right] at (5.8,0.75)
        {\scriptsize one active\\[-1mm]\scriptsize direction};
    
    % Bottom plot: singular values after integration
    \begin{scope}[yshift=-2.5cm]
    
    \draw[->] (0,0) -- (7.2,0) node[right] {index $i$};
    \draw[->] (0,0) -- (0,1.25) node[above] {$\sigma_i$};
    
    % Bars
    \foreach \j in {1,...,6} {
        \draw[fill=blue!35,draw=blue!70] 
            (\j-0.25,0) rectangle (\j+0.25,0.75);
    }
    
    % Local threshold h^2
    \draw[dashed,red!80,line width=1pt] (0.3,0.55) -- (6.8,0.55);
    \node[red!80,right] at (6.8,0.55) {\scriptsize fine tolerance $2h \xi_n$};
    
    % Global tolerance annotation
    \draw[decorate,decoration={brace,amplitude=5pt,mirror}]
        (0.7,-0.25) -- (6.3,-0.25)
        node[midway,below=6pt,align=center]
        {\scriptsize many small singular values\\[-1mm]
         \scriptsize but total norm $\leq \xi_n$};
    
    \node[align=center] at (3.5,1.15)
        {\small $\bar r_n = 0$, but $\mathcal R_{2h\xi_n}$ keeps many modes};
    
    \end{scope}
    
    \end{tikzpicture}
    \caption{Schematic representation of the example in Section \ref{sec: canc smal}. Although $F$ has rank one at every time, successive pulses generate many small singular values. Each of them exceeds the fine recompression tolerance $2h \xi_n = \mathcal O(h^2)$, while their total contribution is below the target accuracy $\xi_n$. Hence all best approximation ranks vanish, whereas the numerical rank grows like $h^{-\beta}$.}
    \label{fig: counter ex 2}
\end{figure}
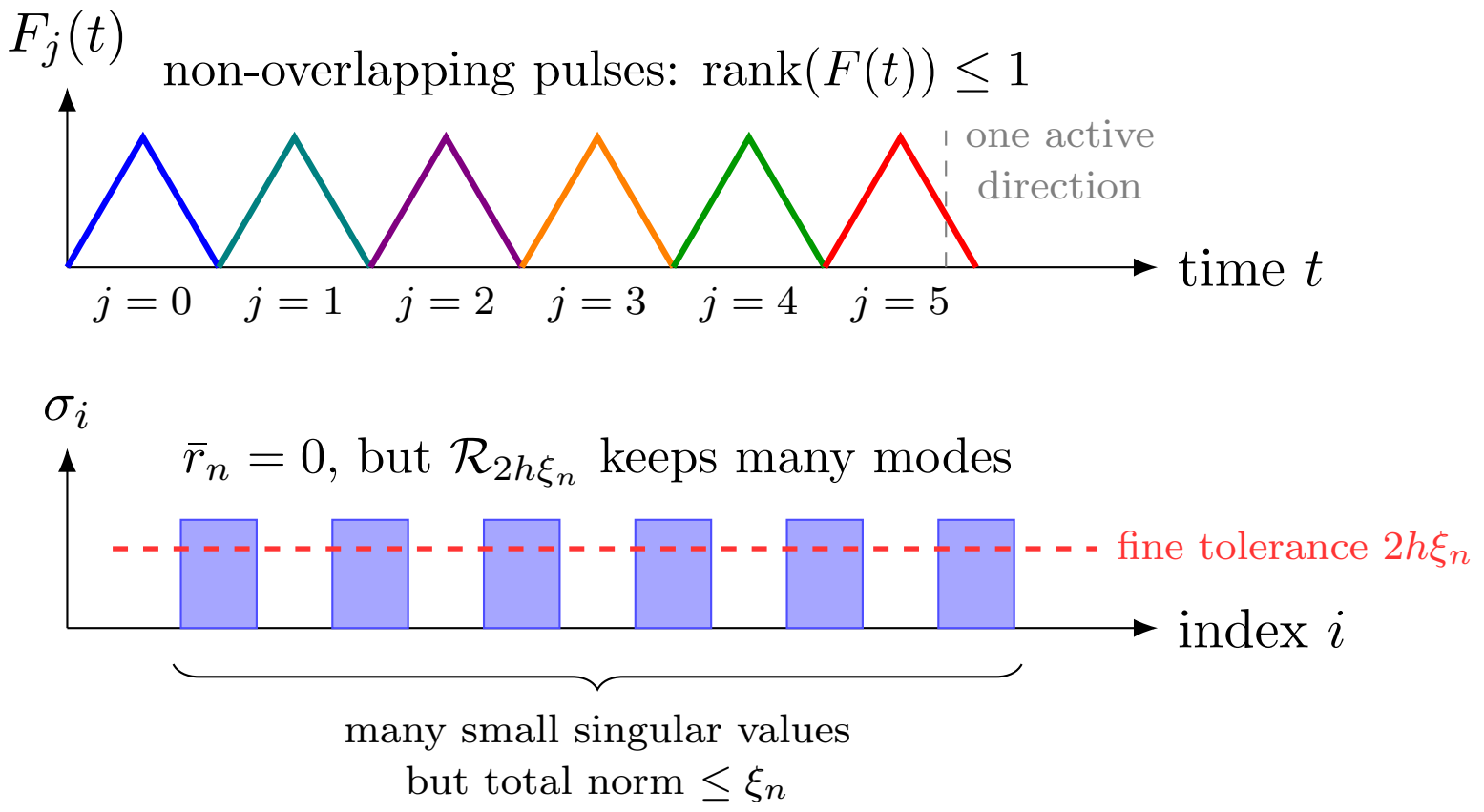

\begin{proof}
    On the time interval $[0,1]$, we have $j_n = 0$ and thus
    \begin{equation*}
        \frac h2 \leq \xi_n = \frac h2 e^{3(\Lambda +1)t_n} \leq K h, \quad \text{for} \quad K = \frac 12 e^{3(\Lambda +1)}.
    \end{equation*}

$(i)$ For $p < \infty$, the assumption on $\beta$ yields $2\alpha - \beta/p > 1$, so that using \eqref{eq: ex 2 bound} we obtain $\| X(t_n)\|_p \leq h/4 \leq \xi_n$ for $h \leq 1$. For $p = \infty$, we have $\| X(t_n) \|_\infty \leq  h^{2 \alpha}/4 \leq \xi_n$. The zero matrix is thus an admissible approximation of the true solution with accuracy $\xi_n$, such that $\bar r_n = 0$. In this setting, it also holds that $R_\Psi = 0$ and $R_0(t, \gamma) \leq \rank \big(F(t)\big) \leq 1$ for all $\gamma \geq 0$. Inserting this into the rank bound of Theorem \ref{thm:main} yields
        \begin{equation*}
            \rank(X_{n+1}) \leq \frac 1h \big( 1+ \max_k R_0(t_k, \xi_k)\big) + R_0(t_n,0) \leq \frac 2h +1, \qquad 0 \leq n < M.
        \end{equation*}
        
$(ii)$ In the given setting, one has $X_1 = 0$ and the algorithm reads
        \begin{equation*}
            X_{n+1} = \mathcal R_{2h \xi_n}(X_n) + hF(t_n), \qquad 1 \leq n \leq M-1,
        \end{equation*}
        since none of these steps is a reset step for $h \leq 1$. By induction, all $X_n$ with $n \leq M$ are diagonal with nonnegative entries, since $F \geq 0$. As these entries correspond to the singular values of $X_n$, the recompression $\mathcal R_\delta$ can only discard entries, each of which is at most $\delta$.

        Next, fix $j \in \{0, \dots, J_h-1 \}$. Since $H(h^{-\alpha}t-j) \geq 1/4$ on the interval
        \[
        I_j = \bigg[h^\alpha\Big(j+\frac 14\Big), h^\alpha\Big(j + \frac 34\Big)\bigg],
        \]
        the $(j+1)$-th entry of $F$ is at least $h^\alpha /4$ on $I_j$. Moreover, since $I_j$ is of length $h^\alpha/2$, there is an index $n_j \geq 1$ with $t_{n_j} \in I_j$ for $h$ sufficiently small. The $(j+1)$-th entry of the approximation $X_{n_j+1}$ is thus at least $h \cdot h^{\alpha}/4 = h^{\alpha+1}/4$. Since $\alpha <1$, it also holds that
        \begin{equation*}
            2h \xi_n \leq 2Kh^2 < \frac{h^{1+ \alpha}}{4}, \qquad n \leq M,
        \end{equation*}
        for sufficiently small $h$, where the first inequality comes from the definition of $K$. Consequently, this entry is never discarded by the fine recompression $\mathcal R_{2 h \xi_n}$ for all $n \leq M$ and since $F \geq 0$, it also does not decrease. Finally, the last pulse ends at $t=J_h h^\alpha$, which implies that for all $n$ with $J_h h^\alpha \leq t_n \leq 1$ the approximation $X_n$ has at least $J_h$ nonzero diagonal entries, see also Figure \ref{fig: counter ex 2}.

        For the last statement in the proposition, let $\beta^\prime < 1$, choose $\beta \in (\beta^\prime,1)$ and $\alpha \in \big( \max\{ \beta, \frac 12 (1 + \beta/p) \},1 \big)$. Then $\beta < \alpha$ and $\beta < p(2\alpha -1)$, which yields $\rank(X_n) \geq \lfloor h^{-\beta} \rfloor$ by $(ii)$. This exceeds $ch^{- \beta^{\prime}}$ for any constant $c > 0$ and $h$ sufficiently small.

\end{proof}

Note that all pulses occur before the first rank reset at $t=1$, which then removes the accumulated entries since $\bar r_M =0$. A possible modification would be to replace the hard recompression by shrinkage-type thresholding, where singular values are not only discarded but also reduced, so that repeated recompressions act like a single recompression with a larger threshold. For instance, the soft thresholding operator $\mathcal S_\delta$, see for example \cite{bachmayr2017iterative}, acting on singular values by $\sigma_k \mapsto \max\{0, \sigma_k-\delta \}$ satisfies $\mathcal S_\delta\circ \mathcal S_\delta = \mathcal S_{2\delta}$, which can help to remove small singular values gradually. Nevertheless, this does not solve the problem in general: if the pulses are shifted towards the end of a reset period, too few recompression steps remain to reduce the newly created singular values below the target accuracy.

\subsection{The recompression may miss cancellations in the true solution}\label{sec: miss canc}

The last example illustrates a limitation of error analyses based on the accumulation of local errors when the low-rank structure of the exact solution is caused by cancellations that occur only at later times.

We therefore consider schemes of the general form
\begin{equation*}
    X_{n+1} = \mathcal T_n(X_n + F_n), \qquad X_0 = 0,
\end{equation*}
where $\mathcal T_n$ is an arbitrary truncation operator and $F_n$ is an approximation of the exact increment $\int_{t_n}^{t_{n+1}} F(s) \, ds$ with the local errors
\begin{equation*}
    \eta_n = \| \mathcal T_n(X_n + F_n) - (X_n + F_n) \|, \qquad \zeta_n = \Big\| F_n - \int_{t_n}^{t_{n+1}} F(s) \, ds \Big\|.
\end{equation*}
This framework now also allows, for example, removing singular values larger than $h^2$ in a single step. Since the right-hand side $F$ does not depend on the solution, the time stepping error satisfies
\begin{equation*}
    X_{n+1} - X(t_{n+1}) = X_n - X(t_n) + \big( \mathcal T_n(X_n + F_n)-(X_n +F_n) \big) + \Big( F_n - \int_{t_n}^{t_{n+1}} F(s) \, ds \Big),
\end{equation*}
and consequently also for all $0 \leq k < \ell \leq N$,
\begin{equation*}
    \big\| (X_\ell - X(t_{\ell})) - (X_k - X(t_{k})) \big\| \leq \sum_{n = k}^{\ell -1}( \eta_n + \zeta_n).
\end{equation*}
In particular, it also holds that $\|X_n - X(t_n) \| \leq \sum_{k < n}(\eta_k + \zeta_k)$.

Following the framework of the previous subsection, we let $\frac 12 < \alpha < 1$, $2 \alpha - 1 < \beta < \alpha$, and define for $R \geq 2$ the forcing term
\begin{equation*}
    F(t) = \sum_{j = 0}^{J_h-1} h^\alpha \big( H(h^{-\alpha}t-j) - H(h^{-\alpha}t-j-2R) \big)E_{j+1}.
\end{equation*}
In other words, the $(j+1)$-th diagonal entry of $F$ has a positive pulse supported on the interval $[h^{\alpha}j,h^{\alpha}(j+1)]$, and a negative one with support $[h^{\alpha}(j+2R),h^{\alpha}(j+2R+1)]$. Since both have the same area, but opposite sign, their total integrals cancel out. At any time $t$, by construction at most one positive pulse and at most one negative pulse are active, hence $\rank(F(t)) \leq 2$. Moreover, $F$ is Lipschitz continuous with constant $\sqrt{2}$ with respect to the Frobenius norm. For $h$ sufficiently small, all pulses occur within $[0,1)$.

As the integrals are canceled out in the exact solution, the $(j+1)$-th component of $X(t)$ can be nonzero only during the time interval $[h^\alpha j, h^\alpha(j+ 2R +1)]$. Thus, for each $t$, at most $2R+2$ diagonal components of $X(t)$ can be nonzero, where the size of each component is bounded by the total area of one positive pulse. Hence, 
\begin{equation*}
    \| X(t) \|_{F} \leq \frac 14 h^{2 \alpha} \sqrt{2R + 2} \leq h
\end{equation*}
for fixed $R$ and sufficiently small $h$ since $\alpha > \frac12$ by assumption, which implies ${\mathcal R_h(X(t))=0}$. In particular, the exact solution has best approximation rank zero at the global accuracy scale $h$.

\begin{figure}
    \centering
    \begin{tikzpicture}[x=1cm,y=0.8cm,>=Latex]
    
    % Time axis
    \draw[->] (0,-0.5) -- (8.5,-0.5) node[right] {$t$};
    
    % Draw several modes
    \foreach \i in {1,...,5} {
    
        \pgfmathsetmacro{\y}{5-\i}
        \pgfmathsetmacro{\x}{0.6*(\i-1)}
    
        % baseline for each mode
        \draw[gray!30] (0,\y) -- (8,\y);
        \node[left] at (0,\y) {\scriptsize mode $\i$};
    
        % positive hat
        \draw[green!60!black,line width=1.1pt]
            (\x,\y) -- (\x+0.35,\y+0.45) -- (\x+0.7,\y);
        \node[green!50!black] at (\x+0.35,\y+0.65) {\scriptsize $+$};
    
        % plateau: component must be stored until cancellation
        \draw[blue!60,line width=1.4pt]
            (\x+0.7,\y) -- (\x+3.2,\y);
    
        % negative hat
        \draw[orange!90!black,line width=1.1pt]
            (\x+3.2,\y) -- (\x+3.55,\y-0.45) -- (\x+3.9,\y);
        \node[orange!90!black] at (\x+3.55,\y-0.65) {\scriptsize $-$};
    }
    
    % Rank budget line
    \draw[dashed,red!80,line width=1pt] (3.15,-0.7) -- (3.15,4.2);
    \node[red!80,align=center,above] at (3.15,4.2)
        {\scriptsize at this time\\[-1mm]\scriptsize many modes are active};
    
    % Braces for active components
    \draw[decorate,decoration={brace,amplitude=5pt,mirror}]
        (6.85,0) -- (6.85,4)
        node[midway,right=7pt,align=left]
        {\scriptsize active components\\[-1mm]
         \scriptsize exceed rank budget $R$};
    
    % Red cross indicating early truncation
    \draw[red!80,line width=1.2pt]
        (2.25,2) +(-0.15,-0.15) -- +(0.15,0.15);
    \draw[red!80,line width=1.2pt]
        (2.25,2) +(-0.15,0.15) -- +(0.15,-0.15);
    
    \node[red!80,align=center] at (1.05,1.4)
        {\scriptsize truncated\\[-1mm]\scriptsize too early};
    
    % Arrow to later cancellation
    \draw[->,red!80,bend left=20]
        (2.35,2.05) to (4.0,2.0);
    \node[red!80,align=center] at (5.6,2.45)
        {\scriptsize later cancellation\\[-1mm]\scriptsize is missed};
    
    \end{tikzpicture}
    \caption{Schematic representation of the example in Section \ref{sec: miss canc}. A positive contribution in one mode is cancelled later by a negative contribution in the same mode. The exact solution remains low-rank because of these cancellations. A scheme with ranks bounded by $R$ has to discard components before their cancelling pulse arrives, which leads to large local errors.}
    \label{fig: counter ex 3}
\end{figure}
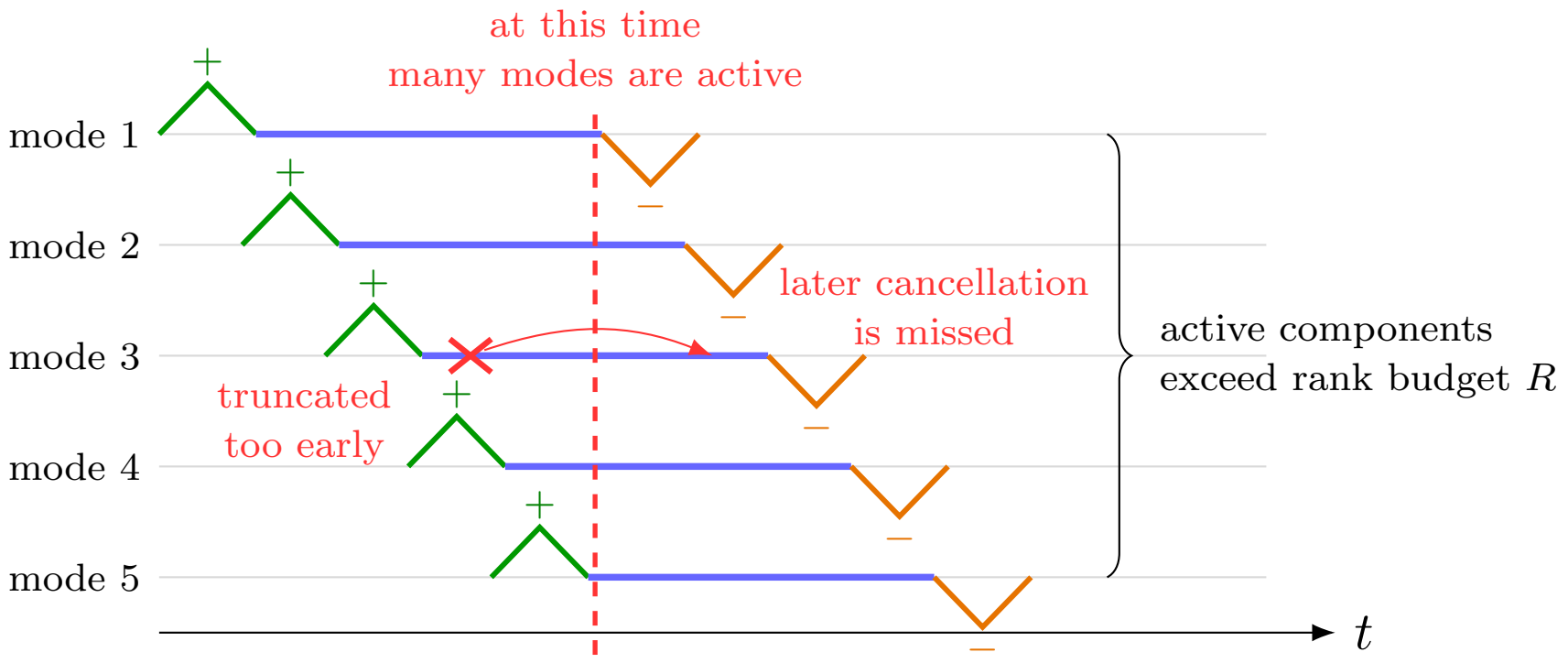

However, a numerical scheme that keeps the rank at most $R$ cannot exploit these cancellations. Each positive pulse is only cancelled after a delay of about $2Rh^\alpha$, during which roughly $2R$ components are waiting for their negative pulse. To keep the rank at most $R$, about half of them must be removed, or never be built up, before their respective cancelling pulse arrives, see Figure \ref{fig: counter ex 3} for $R=2$. Since each positive pulse has size $h^{2\alpha}/4$ and there are $J_h$ of them, this heuristically leads to accumulated truncation and forcing errors of order
\begin{equation*}
    \sum_{n=0}^{\frac 1h -1} (\eta_n + \zeta_n) \gtrsim J_hh^{2\alpha} \approx h^{2\alpha - \beta} \gg h, \quad \text{for} \quad h \rightarrow 0.
\end{equation*}

Hence an error analysis based on the accumulation of local errors cannot certify a scheme of bounded rank in this example. Note that this concerns the analysis rather than the accuracy: the removed contributions may cancel later, but this cannot be detected in the individual time steps.

In summary, Proposition \ref{prop: h sharpness} shows that the factor $1/h$ in Theorem \ref{thm:main} cannot be improved to $h^{-\beta}$ with $\beta < 1$ for Algorithm \ref{alg:step_trunc_reset}. Section \ref{sec: wrong sing} shows that classical step truncation schemes with tolerances at the scale of the local error may lead to ranks of order $1/h$ even for problems that Algorithm \ref{alg:step_trunc_reset} treats with bounded ranks, and Section \ref{sec: miss canc} illustrates that error analyses based on the accumulation of local errors cannot certify $h$-independent rank bounds in general.

\section{Numerical experiment: time-varying explicit DLRA decomposition}\label{sec: num exp}

In this section, we visualize the performance of the proposed integrator on a numerical experiment and compare it to the rank-adaptive BUG integrator. All implementations were carried out using Julia 1.13.3 as well as a fixed tolerance $\xi$ specified below instead of defining it as in \eqref{eq: thm 1 opt ranks} for simplicity. We restrict our considerations to the Frobenius norm case $p=2$, for square matrices of size $d_1=d_2=d$.

All low-rank operations in the proposed step-truncation schemes are carried out in factored form: given $Y_1 = U_1 S_1 V_1^T$ and $Y_2 = U_2 S_2 V_2^T$ of ranks $r_1$ and $r_2$, their sum is given by a concatenation of the respective factors, that is,
\begin{equation}\label{eq: Y_1 + Y_2}
    Y_1 + Y_2 = \begin{bmatrix} U_1 & U_2 \end{bmatrix}
    \begin{bmatrix} S_1 & 0 \\ 0 & S_2 \end{bmatrix}
    \begin{bmatrix} V_1 & V_2 \end{bmatrix}^T,
\end{equation}
which yields a representation of rank $r_1 + r_2$ with non-orthonormal factors. However, in general this representation is redundant and we thus include an additional recompression to a fine tolerance after each addition of two low-rank matrices.

Sums of low-rank matrices are recompressed without forming the full $d\times d$ matrix in \eqref{eq: Y_1 + Y_2}. To this end, we compute reduced QR decompositions
\[
\begin{bmatrix} U_1 & U_2 \end{bmatrix} = Q_U R_U, \qquad \begin{bmatrix} V_1 & V_2 \end{bmatrix} = Q_V R_V,
\]
and the singular value decomposition of the small core matrix
\[
C = R_U \begin{bmatrix} S_1 & 0 \\ 0 & S_2 \end{bmatrix} R_V^\top = \hat U \hat\Sigma \hat V^\top \in \R^{(r_1+r_2)\times(r_1+r_2)}.
\]
Since $Q_U\hat U$ and $Q_V\hat V$ have orthonormal columns, $Y_1+Y_2 = (Q_U\hat U)\,\hat\Sigma\,(Q_V\hat V)^\top$ is a singular value decomposition, which we truncate according to \eqref{eq: R_delta} with $p=2$, that is, we compute $\mathcal R_\delta(Y_1+Y_2)$. This requires only $\mathcal O\big(d(r_1+r_2)^2\big)$ operations. For these intermediate recompressions, we use the tolerance $\delta = \xi h/2$ in the following.

As a numerical example, we consider a manufactured matrix ODE of the standard type considered in the dynamical low-rank literature: we first specify an exact solution of the form $X(t) = U(t) \Sigma(t) V(t)^T \in \mathbb R^{d \times d}$ with analytically prescribed orthonormal bases and time-dependent singular values and then choose the right-hand side accordingly.

More specifically, we let the basis matrices evolve through rotations generated by sparse, block-structured skew-symmetric matrices $\Omega_U$ and $\Omega_V$, that is,
\begin{equation*}
    \dot{U}(t) = \Omega_U U(t), \quad \dot{V}(t) = \Omega_V V(t),
\end{equation*}
where
\begin{equation*}
    \Omega_U = \sum_{j=1}^\ell \omega_{U,j}(e_{\ell+j}e_j^T - e_j e_{\ell+j}^T), \quad \Omega_V = \sum_{j=1}^\ell \omega_{V,j}(e_{\ell+j}e_j^T - e_j e_{\ell+j}^T),
\end{equation*}
for some angular velocity vectors $\omega_U, \omega_V \in \mathbb R^\ell$ with $2\ell \leq d$. Here, $e_1, \dots, e_d \in \mathbb R^d$ denote the canonical basis vectors. We then have an explicit representation
\begin{align*}
    u_j(t) &= e^{t \Omega_U}e_j = \cos(\omega_{U,j}t)e_j + \sin(\omega_{U,j}t)e_{\ell+j}, \\
    v_j(t) &= e^{t \Omega_V}e_j = \cos(\omega_{V,j}t)e_j + \sin(\omega_{V,j}t)e_{\ell+j},
\end{align*}
for $1 \leq j \leq \ell$, of the time-dependent columns of the orthonormal matrices $U(t)$ and $V(t)$.

We further prescribe the time-varying singular values stored in $\Sigma(t) = \mathrm{diag}(\sigma_1(t), \dots, \sigma_\ell(t))$ through
\begin{equation*}
    \sigma_j(t) = \exp(-\alpha(t)(j-1)), \qquad 1 \leq j \leq \ell,
\end{equation*}
for
\begin{equation}\label{eq: alpha_t num ex1}
    \alpha(t) = \alpha_{\max} - ( \alpha_{\max} - \alpha_{\min} ) \sin^2 \Big(\frac{2 \pi t}{T} \Big), \qquad \alpha_{\min} = 0.25, \quad \alpha_{\max} = 0.8.
\end{equation}
Using this construction, we obtain the ODE
\begin{equation*}
    \dot{X}(t) = F(t,X) = \Omega_UX - X \Omega_V + G(t), \quad G(t) = U(t) \dot{\Sigma}(t) V(t)^T,
\end{equation*}
where $\dot{\sigma}_j = (j-1)(\alpha_{\max}-\alpha_{\min})\frac{2\pi}{T}\sin \Big( \frac{4 \pi t}{T} \Big) \sigma_j(t)$ for $1 \leq j \leq \ell$, with the initial data
\begin{equation*}
    X(0) = \mathrm{diag}(1,e^{-\alpha_{\max}},e^{-2\alpha_{\max}},\dots,e^{-(\ell-1)\alpha_{\max}},0, \dots,0).
\end{equation*}
We further set $d = 200$, $\ell=50$, $h = 0.01$, $T=10$ and choose the $\ell$ entries of $\omega_U$ and $\omega_V$ to be distributed equidistantly in $[0.3,1.3]$ and $[0.4,1.1]$, respectively. In this setting, the smallest singular value is of size $e^{-0.25 \cdot 49} \approx 5 \cdot 10^{-6}$. The initial data is truncated to the corresponding approximation tolerance $\xi \in \{ 10^{-4}, 10^{-5} \}$. 

\begin{figure}[!h]
   	\centering
   	\includegraphics[width = \linewidth]{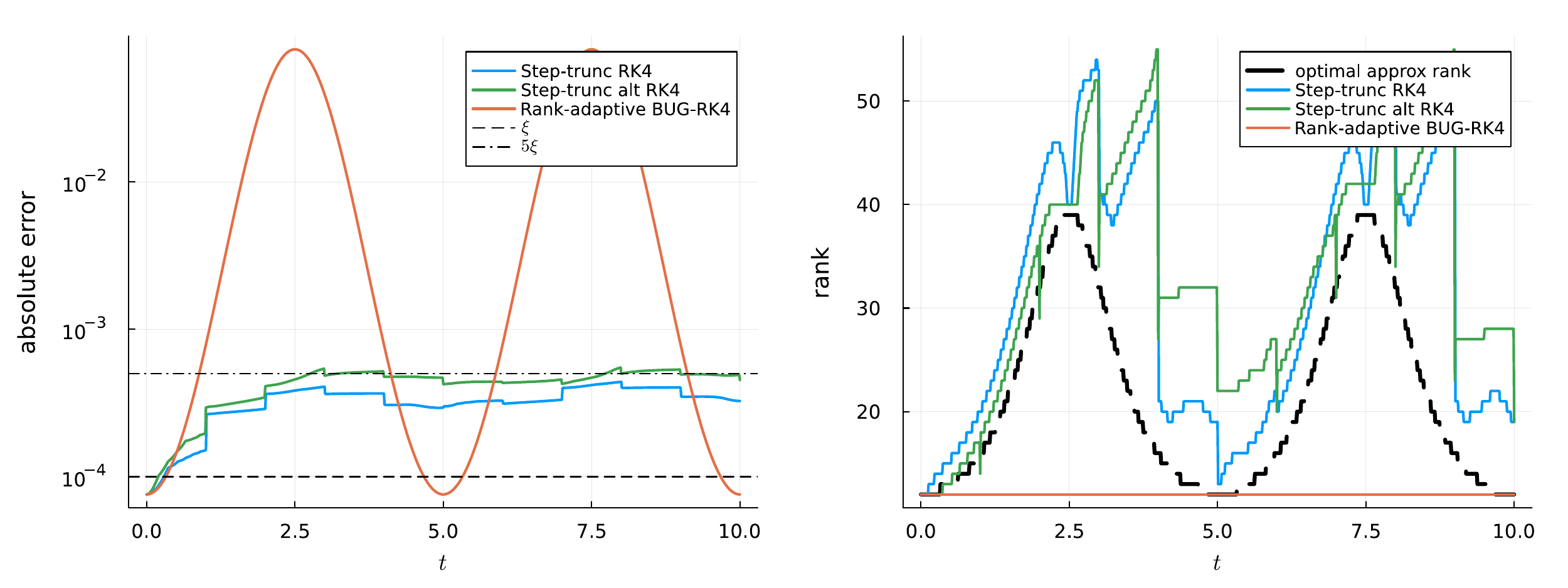}

    \includegraphics[width = \linewidth]{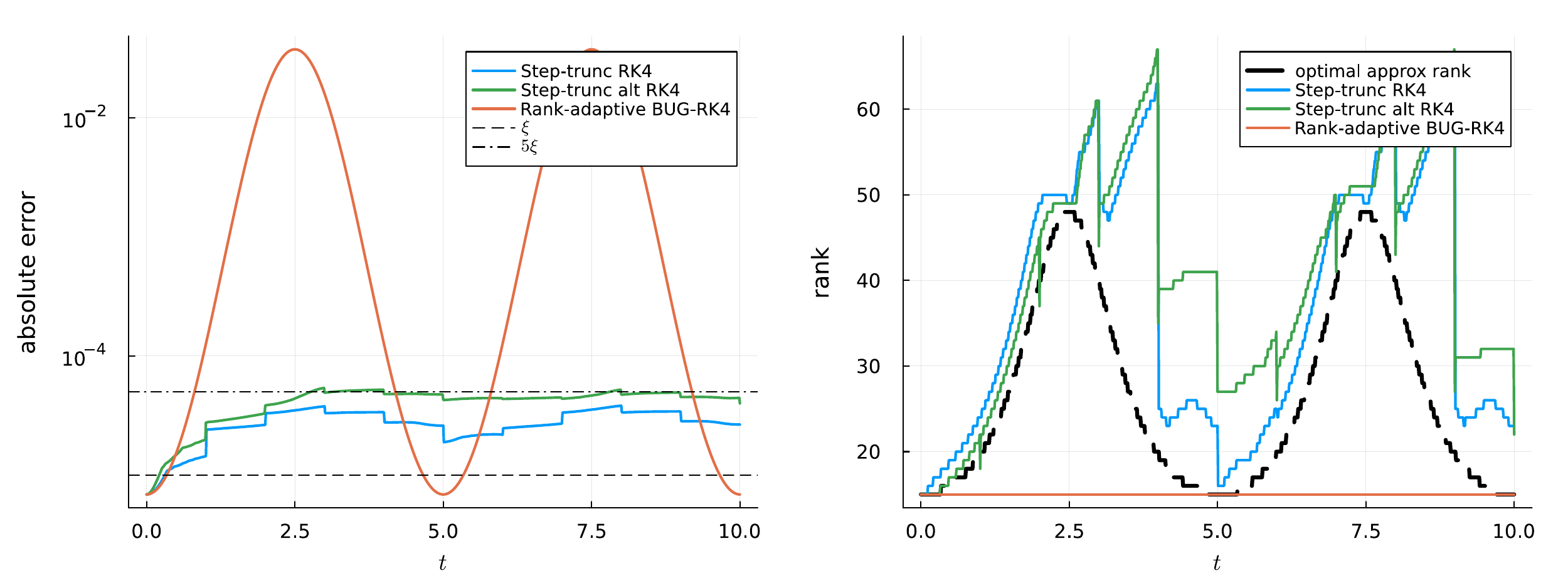}
   		\caption{Time-varying explicit DLRA decomposition: comparison of the absolute error in Frobenius norm with respect to the true solution using different time steppers (left) and the corresponding evolution of approximation ranks compared to the optimal ones (right) using a tolerance of $\xi = 10^{-4}$ (top) or $\xi = 10^{-5}$ (bottom).}
    \label{fig: numerics1}
\end{figure}

Figure \ref{fig: numerics1} visualizes the evolution of the absolute error in Frobenius norm with respect to the true solution as well as the time evolution of the approximation ranks for the two proposed algorithms `Step-trunc RK4' (Algorithm \ref{alg:step_trunc_reset}) and `Step-trunc alt RK4' (Remark~\ref{rem: alt thresh}) using the explicit RK4 method as the base method. The methods are compared to the rank-adaptive BUG method, for which we use the same truncated initial data and step size $h$, solve the K-, L- and S-substeps with the explicit RK4 method, and choose the truncation tolerance as $\vartheta = \xi h$.

\begin{figure}
\centering
\begin{minipage}{.48\linewidth}
  \centering
  \includegraphics[width=\linewidth]{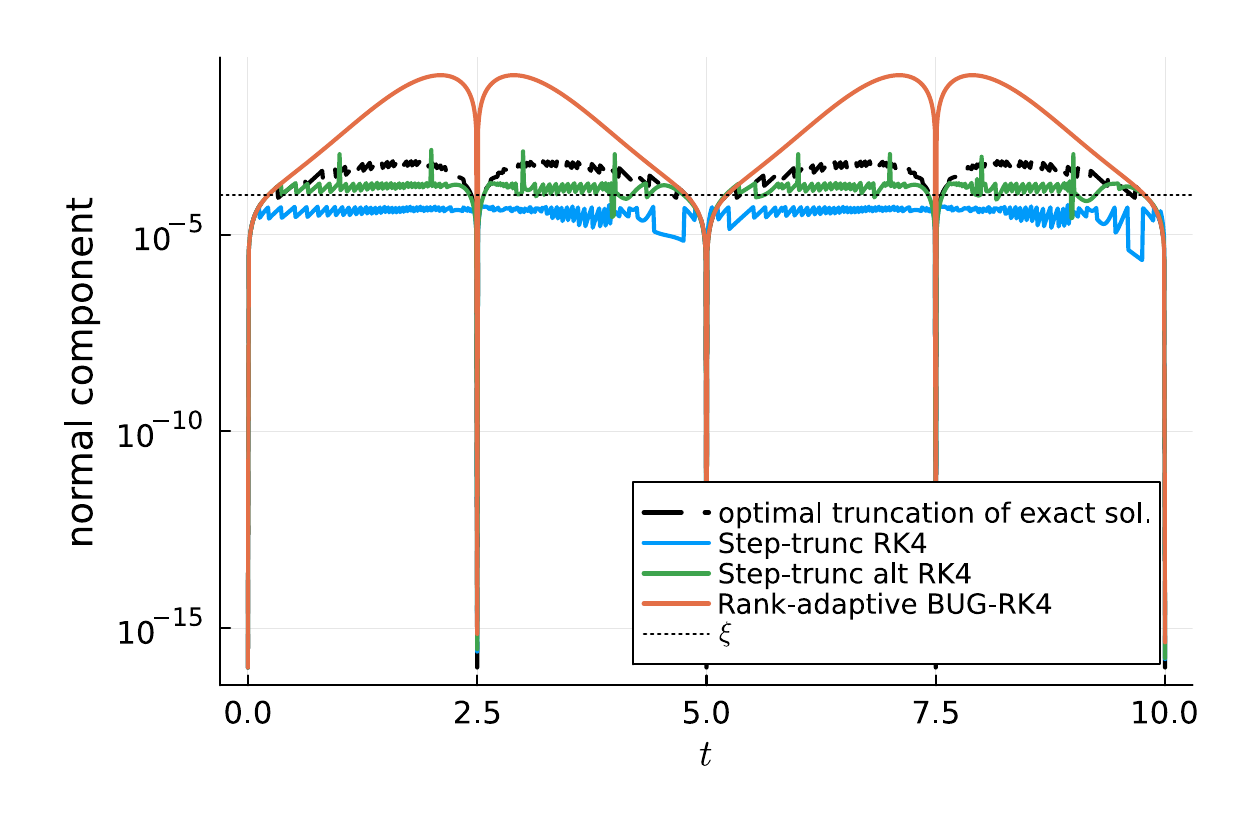}
\end{minipage}%
\begin{minipage}{.48\linewidth}
  \centering
  \includegraphics[width=\linewidth]{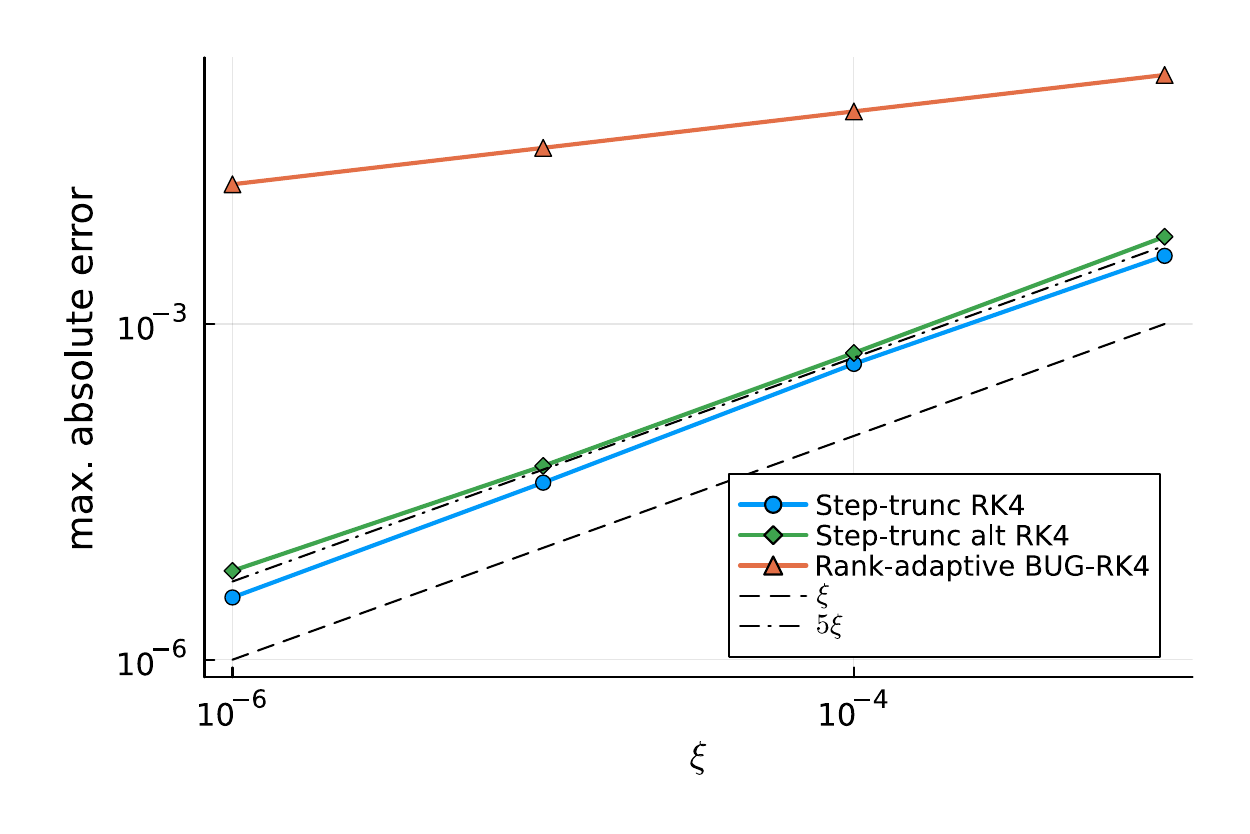}
\end{minipage}
\caption{Visualization of the size of the normal component in Frobenius norm over time using $\xi = 10^{-4}$ (left) and the maximal absolute error over time as a function of $\xi$ (right).}
\label{fig: numerics 2}
\end{figure}

While both proposed algorithms achieve an accuracy in $\mathcal{O}(\xi)$, the absolute error associated with the rank-adaptive BUG integrator appears to be multiple orders of magnitude higher. The ranks of the rank-adaptive BUG integrator remain constant at the ranks $r_0 = 12$ and $r_0 = 15$ of the truncated initial data, and its error exceeds the tolerance by about three orders of magnitude. This behavior can be explained exactly: for $Y = \sum_{j \leq r} \sigma_j u_j v_j^T$, the rotation part of the right-hand side $F$ lies in the tangent space $\mathcal T_Y \mathcal M_r$, so that the
normal component is given entirely by the source term, that is,
\begin{equation*}
    (I - P(Y))F(t,Y) = \sum_{j = r+1}^\ell u_j(t)\, \dot{\sigma}_j(t)\, v_j(t)^T,
\end{equation*}
where $P(Y)$ denotes the orthogonal projection onto the tangent space at $Y$ of the manifold of rank-$r$ matrices. Let $\mathcal W_{r_0} = \Span\{e_j, e_{\ell+j} : 1 \leq j \leq r_0\}$. Since the rotation planes $\Span\{e_j, e_{\ell+j}\}$ are mutually orthogonal, we have $u_j(t), v_j(t) \in \mathcal W_{r_0}$ for $j \leq r_0$ and $u_j(t), v_j(t) \perp \mathcal W_{r_0}$ for $j > r_0$, and $\Omega_U$, $\Omega_V$ leave $\mathcal W_{r_0}$ invariant. Now assume that the columns of the basis matrices $U_n$, $V_n$ lie in $\mathcal W_{r_0}$, which holds for the truncated initial data. Then $G(t) V_n$ and $G(t)^\top U_n$ only involve $u_j$, $v_j$ with $j \leq r_0$, so that the K- and L-steps, and hence the augmented bases, remain in $\mathcal W_{r_0}$. The Galerkin step only updates the core matrix, and the subsequent truncation selects a subspace of the augmented bases, so that the columns of $U_{n+1}$, $V_{n+1}$ again lie in $\mathcal W_{r_0}$. By induction, the column and row spaces of all BUG iterates $Y_n$ are contained in $\mathcal W_{r_0}$, which yields the lower bound
\begin{equation*}
    \| Y_n - X(t_n) \|_F \geq \Big( \sum_{j = r_0+1}^{\ell} \sigma_j(t_n)^2 \Big)^{1/2},
\end{equation*}
independently of the step size and the truncation tolerance. The augmented
directions thus only capture the rotation of the singular vectors already
present and are partially removed again by the truncation. This does not contradict the error bound of the rank-adaptive BUG integrator, which assumes the normal component $\|(I - P(Y))F(t,Y)\|$ to be small uniformly in time, see for example \cite[Equation (3.12)]{ceruti2024parallel}. In the present example it exceeds $\xi$ in both increasing and decreasing rank phases, except near the turning points $t \in \{0,2.5,5.0,7.5,10 \}$, see the left-hand side of Figure~\ref{fig: numerics 2}. Note that the normal
component vanishes at the turning points of $\alpha$ in \eqref{eq: alpha_t num ex1}, where $G(t) = 0$: it measures the rate at which the exact solution leaves the approximation space rather than the accumulated deviation.

Comparing the two proposed step truncation methods, the variant given in Algorithm~\ref{alg:step_trunc_reset} can be observed to be slightly more accurate than the one discussed in Remark \ref{rem: alt thresh}. Concerning the ranks of the produced iterates, those of the main scheme are lower after resets and in decreasing phases. The alternative is slightly lower when the ranks increase. Additionally, the first variant is almost able to reproduce the lowest optimal approximation ranks at $t=5.0$ and we observe a near-periodic reproduction of the approximation ranks over time for both algorithms. Note that here however the reproduction of the decreasing rank behavior is possible only because of the periodic rank resets every 100 time steps, which was confirmed numerically, underlining the importance of this reset step. The right-hand side of Figure \ref{fig: numerics 2} shows the dependence of the actual maximal absolute error in Frobenius norm on the chosen tolerance $\xi$: indeed, the plot indicates that the error remains well within $\mathcal O(\xi)$ for the tested tolerances.

\section{Conclusion and Perspectives}

In this work, we proposed a rank-adaptive integrator for ODEs, based on step-truncation procedures of an explicit Runge-Kutta method, requiring no prior knowledge of the exact solution. This scheme exploits two recompression operations based on different tolerances: the former acts at the level of the increment function of the Runge-Kutta method, whereas the latter acts on the propagated iterate. This strategy, together with the presented choice of truncation tolerances, permits us to reach near-optimal ranks, up to a factor of $1/h$, at any integration time, while keeping the same order of accuracy of the underlying Runge-Kutta procedure. Moreover, the provided counterexamples illustrate the optimality of the obtained rank bound with respect to the dependence on the inverse of the time step size $h$. The numerical results illustrate  efficient and accurate performance in practice, while also providing a comparison with the rank-adaptive BUG integrator.

Natural extensions of this work concern, for instance, rank-adaptive integration techniques for stochastic differential equations. Moreover, an extension of the results to higher-dimensional tensor valued ordinary differential equations as well as a further investigation of the exact placement of possible additional recompressions are already in preparation.

\section{Acknowledgements}

The authors thank Markus Bachmayr for helpful discussions and comments on a draft of this manuscript, and Youssef Marzouk for hosting M.D. and F.Z. in his research group, where the idea for this collaboration arose. P.S. acknowledges funding by the European Union (ERC, COCOA, 101170147). Views and opinions expressed are however those of the authors only and do not necessarily reflect those of the European Union or the European Research Council. Neither the European Union nor the granting authority can be held responsible for them. F.Z. acknowledges the support of Swiss National Science Foundation under the Doc.Mobility funding linked to Project n. 200518 “Dynamical low rank
methods for uncertainty quantification and data assimilation”. F.Z. also thanks Fabio Nobile for travel support. M.D. and P.S. were also supported by the National Science Foundation under Grant No. DMS-1929284 and were in residence at the Institute for Computational and Experimental Research in Mathematics in Providence, RI, during the ``Stochastic and Randomized Algorithms in Scientific Computing: Foundations and Applications" semester program.
We acknowledge the use of the \emph{LowRankTensors} julia package (see \url{https://gitlab.git.nrw/rwth-igpm/igpmregistry/-/tree/main/L/LowRankTensors}) developed by Markus Bachmayr's group. ChatGPT 5.5 and Claude Opus 5.5 were used for literature research, proofreading, and coding assistance in Julia and TikZ.

\printbibliography

\end{document}